\documentclass[11pt, reqno]{amsart}
\usepackage{epsfig}

\usepackage{amstext}
\usepackage{xspace}
\usepackage{amsfonts}
\usepackage{amsmath}
\usepackage{amssymb}
\usepackage{amstext}
\usepackage{amsthm}       
\usepackage{xspace}

\newcommand{\C}{\mathbb C}
\newcommand{\EE}{\mathbb E}
\newcommand{\FF}{\mathbb F}
\newcommand{\G}{\mathbb G}
\newcommand{\R}{\mathbb R}

\newcommand{\N}{\mathbb N}

\newcommand{\T}{\mathbb T}
\newcommand{\Z}{\mathbb Z}
\newcommand{\restr}{\hbox{\LARGE $\llcorner$}}

\newcommand{\comp}{\mbox{\scriptsize  $\circ$}}
\newcommand{\eps}{\varepsilon}

\newcommand{\ucv}{\rightrightarrows}

\newcommand{\A}{\mathcal{A}}
\newcommand{\hh}{\mathcal{H}}

\newcommand{\e}{\textrm{\rm e}}
\newcommand{\Eq}{\mathcal{E}}
\newcommand{\F}{\mathcal{F}}

\newcommand{\M}{\mathcal{M}}
\renewcommand{\S}{\mathcal{S}}

\newcommand{\leb}{\mathcal{L}}

\newcommand{\PP}{\mathbb{P}}
\newcommand{\tagliato}{$\kern-5.5 mm -$}
\newcommand{\tagliat}{$\kern-6.5 mm -$}

\newcommand{\tagli}{\mbox{\tagliat}}
\newcommand{\cchi}{\mbox{\large $\chi$}}
\newcommand{\abra}[1]{(\ref{#1})}

\newcommand{\D}[1]{\mbox{\rm #1}}
\newcommand{\dd}{\D{d}}

\newtheorem{teorema}{Theorem}[section]
\newtheorem{prop}[teorema]{Proposition}
\newtheorem{lemma}[teorema]{Lemma}
\newtheorem{definition}[teorema]{Definition}
\newtheorem{cor}[teorema]{Corollary}

\newtheorem{guess}[teorema]{Remark}
\newtheorem{example}[teorema]{Example}

\newenvironment{dimo}{{\bf\noindent Proof.}}{\qed}
\newenvironment{oss}{\begin{guess} \begin{rm}}{\end{rm} \end{guess}}
\newenvironment{definizione}{\begin{definition} \begin{rm}}{\end{rm}
\end{definition}}

\begin{document}

\title{ Weak KAM Theory topics\\
in the stationary ergodic setting}
\author{Andrea Davini \and Antonio Siconolfi}
\address{Dip. di Matematica, Universit\`a di Roma ``La Sapienza",
P.le Aldo Moro 2, 00185 Roma, Italy}
\email{davini@mat.uniroma1.it,\quad siconolf@mat.uniroma1.it }

\begin{abstract}
We perform a qualitative analysis of the critical equation
associated with a stationary ergodic Hamiltonian through a
stochastic version  of the metric method, where the notion of
closed random stationary set, issued from stochastic geometry,
plays a major role. Our purpose is to give an appropriate notion
of random Aubry set, to single out characterizing conditions for
the existence of exact or approximate correctors, and write down
representation formulae for them. For the last task,  we make
 use  of  a Lax--type formula,  adapted to the
stochastic environment. This  material  can be regarded as a first
step of a long--term project to develop a random analog of Weak
KAM Theory, generalizing what done in the periodic case or, more
generally, when the underlying space is a compact manifold.\par
\end{abstract}
\maketitle

\begin{section}{Introduction}
For a given a probability space $\Omega$,  on which $\R^N$ acts
ergodically, we consider the family of Hamilton--Jacobi equations
\begin{equation*}\label{eq intro HJa}
H(x,Dv,\omega)=a\qquad\hbox{in $\R^N$,}
\end{equation*}
where $a$ varies in $\R$,  and $H$ is a continuous Hamiltonian,
convex and superlinear in the momentum variable, and stationary
with respect to the action of $\R^N$.  As it is well known, this
framework includes the periodic \cite{LPV}, quasi--periodic
\cite{Ar98} and almost--periodic cases \cite{Is} as particular
instances.

A {\em stationary critical value}, denoted by $c$, can be defined
in this setting as the minimal value $a$ for which the above
equation possesses {\em admissible} subsolutions, that is
Lipschitz random functions that have stationary gradient with mean
0 and that are almost sure subsolutions  either in the viscosity
sense or, equivalently, almost everywhere in $\R^N$. The condition
on the gradient implies almost sure sublinear growth at infinity,
see \cite{DavSic08, DS1-08}. The stationary critical value is in
general distinct  from the {\em free critical value} $c_f$, i.e.
the minimal value $a$ for which the above equation admits
subsolutions, without any further qualification. More precisely
$c_f(\omega)$ is a random variable, almost surely constant because
of the ergodicity assumption. Clearly $c \geq c_f$.

The relevance of the stationary critical value $c$ relies on the
fact that it is the only level of $H$ for which the corresponding
critical equation can have admissible exact or approximate
solutions, also named {\em exact} and {\em approximate correctors}
for the role they play in associated homogenization problems, see
Section \ref{sez HJ} for precise definitions.

The aim of  the  paper is  to perform a qualitative study of the
critical equation, in any space dimensions,  through the metric
approach,   by developing the ideas of \cite{DavSic08, DS1-08}.
The adaptation of this pattern to the stationary ergodic setting
requires   the use of some tools from random set theory, the
leading idea being that the stationary ergodic structure of the
Hamiltonian induces a stochastic geometry in the space of the
state variable, where the fundamental entities are the closed
random stationary sets which, somehow, play the same role as the
points in the deterministic case.

More specifically, our purpose  is to give an appropriate notion
of random Aubry set, to single out characterizing conditions for
the existence of exact or approximate correctors, and write down
representation formulae for them. This material can be regarded as
a first step of a long--term project to develop a random analog of
Weak KAM Theory, generalizing what done in the periodic case or,
more generally, when the underlying space is a compact manifold,
see \cite{FaSic03}.\par

We recall that the  random version of the metric method has
allowed to completely clarify  the setup  in the one--dimensional
case \cite{DavSic08}, where it has been proved the existence of
approximate or exact correctors via Lax representation formulae,
depending on whether $0$ belongs or not to the interior of the
flat part of the effective Hamiltonian obtained via homogenization
\cite{ReTa00, Souga99}. This permits, among other things,  to
carry out the homogenization procedures  through Evans' perturbed
test function method.

Even if in the  multidimensional analysis \cite{DS1-08} many
analogies with the one--dimen- sional setting appear, the topic is
definitely more involved, due to the increased degrees of freedom,
so that the picture is far from being complete. In particular the
issue of the existence of approximate correctors is a relevant
open problem, see Section \ref{sez open}.

Our  investigation   can be briefly described as follows. We
associate to the critical equation a Finsler--type random
semidistance $S$ on $\R^N$, and we consider the family of
fundamental (critical) admissible subsolutions obtained via the
Lax formula
\begin{equation}\label{intro Lax}
\inf\{g(y,\omega)+S(y,x,\omega)\,:\,y\in C(\omega)\,\},
\end{equation}
where $C(\omega)$ is a closed random stationary set and $g$ is an admissible
critical subsolution.

We first address our attention to detect characterizing conditions on $g$
and $C(\omega)$ under which the above formula defines an exact corrector.
  In case $c=c_f$, this holds true   if  $C(\omega)\subset \A_f(\omega)$ almost surely, where $\A_f(\omega)$ is the classical Aubry set,
  made up, as in the deterministic case, by points around which some degeneracy of $S$ takes place.
  It can be defined  through conditions  on cycles, see Section \ref{sez HJ}.
    If instead  $c > c_f$ or $c=c_f$ and $C(\omega)\cap \A_f(\omega) = \emptyset$  a.s., we find that formula \eqref{intro Lax}
    gives a solution if and only
if any point $y_0$ in $C(\omega)$ is connected with the
``infinity'' through a curve along which  $g(\cdot,\omega)$ is
equal to $g(\cdot,\omega)+S(\cdot,y_0,\omega)$. This, in turn,
implies that the asymptotic norm associated to $S$ is degenerate.

The subsequent step is to use this information to propose a
suitable notion of {random Aubry set} and to explore its
properties. Our choice, in analogy with the periodic setting, is
to define the {\em random Aubry set} $\A$ as the maximal
stationary closed random set that plugged into \eqref{intro Lax}
in place of $C$ defines a corrector for any choice of the
admissible subsolution $g$. We find that if $c=c_f$  then
$\A_f(\omega)\subset \A(\omega)$ a.s. and if, in addition, no metric degeneracy occurs
at infinity or, in other term, the stable norm associated with $S$
is strictly positive in any direction, then $\A(\omega)$ and $\A_f(\omega)$
almost surely coincide.

Further we prove, generalizing  a property holding in the
deterministic case,   the existence of an admissible critical
subsolution $\overline v$ which is {\em weakly strict} in
$\R^N\setminus\A(\omega)$, i.e. almost surely satisfying
\[
\overline v(x,\omega)-\overline v(y,\omega)<S(y,x,\omega)
\qquad\hbox{for every $x,\,y\in\R^N\setminus\A(\omega)$\quad with
$x\not=y$.}
\]
We  are also able to extend to the stationary ergodic case  some
dynamical properties of the Aubry set. More precisely,  we show
that the random Aubry set is almost surely foliated by curves
defined in $\R$ along which any critical subsolution agrees with
the semidistance $S$, up to additive constants. These curves turn
out to be global minimizers for the action of the Lagrangian in
duality with $H$, and,  when $H$ is regular enough, they are in
addition integral curves of the Hamiltonian flow.

The results  on the random Aubry set  are obtained under the
crucial hypothesis  that $\Omega$ is {\em separable} from the
measure theoretic viewpoint, meaning that $L^2(\Omega)$ is
separable. This assumption, while standard in the probabilistic
literature, would exclude here the almost--periodic case.
Following the usual approach, in fact, an almost--periodic
function can be seen as the restriction on $\R^N$ of a continuous
map defined on $\G^N$, the {\em Bohr compactification} of $\R^N$.
The associated normalized Haar measure is a probability measure
which is ergodic with respect to the action of $\R^N$. This allows
to include the almost--periodic case within the stationary ergodic
framework, but the problem is that $\G^N$ is non--separable, see
\cite{AmFrid} for similar issues.

Thus we have to resort to a different construction, exposed in the
Appendix, that we believe of independent interest. We basically
exploit  that any almost--periodic function on $\R^N$ is the
uniform limit of a sequence of quasi--periodic functions, which,
in turn,  can be seen as  specific realizations of  stationary
ergodic maps defined on  $k$--dimensional tori, with $k$ suitably
chosen.  \par

 By properly defining the objects we work with, we
obtain  that a given almost--periodic Hamiltonian can be seen as a
specific realization of a stationary ergodic one, with $\Omega$
equal to a  {\em countable} product of finite dimensional tori.
The latter, endowed with the product distance, is a compact metric
space, thus separable both from a topological and a
measure--theoretic viewpoint. Some attention must be paid in the
previous construction in order to preserve the ergodicity of the
action
of $\R^N$ on $\Omega$.\\

\indent The paper is organized as follows: in Section \ref{sez
preliminari} we fix notations and expose some preliminary
material, in particular we present definitions and properties of
stationary closed random sets and random functions that are
relevant  for our analysis. Section \ref{sez HJ} is focused on
stochastic Hamilton--Jacobi equations,  we introduce the metric
tools we will need, and  recall some basic facts about
Aubry--Mather theory in the deterministic setting. Section
\ref{sez Lax} is devoted to Lax formulae in the stationary ergodic
setting, in particular to derive characterizing conditions on the
source set and on the trace under which the corresponding Lax
formula defines an exact corrector. In Section \ref{sez Aubry} we
define the random Aubry set and  study its properties. In Section
\ref{sez open} we discuss some questions left open by our study.
The Appendix contains the construction outlined above, addressed
to include  the almost--periodic case in our framework.\medskip
\end{section}

\smallskip\indent{\textsc{Acknowledgements. $-$}}
The first author has been supported for this research by the
European Commission through a Marie Curie Intra--European
Fellowship, Sixth Framework Program (Contract
MEIF-CT-2006-040267). He wishes to thank Albert Fathi for many
interesting discussions and suggestions.\bigskip
\par

\begin{section}{Preliminaries}\label{sez preliminari}

We write below a list of symbols used throughout this paper.
\[
\begin{array}{ll}
N & \hbox{an integer number}\\
B_R(x_0) &\hbox{the closed ball in $\R^N$
centered at $x_0$ of radius $R$}\\
B_R & \hbox{the closed ball in $\R^k$
centered at $0$ of radius $R$}\\
\langle\,\cdot\;, \cdot\,\rangle & \hbox{the scalar product in
$\R^N$} \\
|\cdot| & \hbox{the Euclidean norm in $\R^N$}\\
\R_+& \hbox{the set of nonnegative real numbers}\\
\mathcal{B}(\R^k) & \hbox{the $\sigma$--algebra of Borel subsets
of $\R^k$}\\
\cchi_E & \hbox{the characteristic function of the set $E$}\\
\end{array}
\]

\vspace{1ex} Given a subset $U$ of $\R^N$, we denote by $\overline
U$ its closure. We furthermore say that $U$ is {\em compactly
contained} in a subset $V$ of $\R^N$ if $\overline U$ is compact
and contained in $V$. If $E$ is a Lebesgue measurable subset of
$\R^N$, we denote by $|E|$ its $N$--dimensional Lebesgue measure,
and qualify $E$ as {\em negligible} whenever $|E|=0$. We say that
a property holds {\em almost everywhere} ($a.e.$ for short) on
$\R^N$ if it holds up to a negligible set. We will write
$\varphi_n\ucv\varphi$ on $\R^N$ to mean that the sequence of
functions $(\varphi_n)_n$ uniformly converges to $\varphi$ on
compact subsets of $\R^N$.\par

With the term {\em curve}, without any further specification, we
refer to a Lipschitz--continuous function from some given interval
$[a,b]$ to $\R^N$. The space of all such curves is denoted by
$\D{Lip}([a,b],\R^N)$, while $\D{Lip}_{x,y}([a,b],\R^N)$ stands
for the family of curves $\gamma$  joining  $x$ to $y$, i.e. such
that $\gamma (a)=x$ and $\gamma (b)=y$, for any fixed $x$, $y$ in
$\R^N$.
%
%
The Euclidean length of a curve
$\gamma$ is denoted by $\hh^1(\gamma)$.\par


\smallskip

Throughout the paper, $(\Omega,\F, \PP)$ will denote a {\em
separable probability space},  where $\PP$ is the probability
measure and $\F$ the $\sigma$--algebra of $\PP$--measurable sets.
Here separable is understood in the measure theoretic sense,
meaning that the Hilbert space $L^2(\Omega)$ is separable, cf.
\cite{Tsi01} also for other equivalent definitions. A property
will be said to hold {\em almost surely} ($a.s.$ for short) on
$\Omega$ if it holds up to a subset of probability 0. We will
indicate by $L^p(\Omega)$, $p\geq 1$, the usual Lebesgue space on
$\Omega$ with respect to $\PP$. If $f\in L^1(\Omega)$, we write
$\EE(f)$ for the mean of $f$ on $\Omega$, i.e. the quantity
$\int_\Omega f(\omega)\,\dd \PP(\omega)$.

We qualify  as {\em measurable} a map from $\Omega$ to itself,  or
 to a topological space $\M$ with Borel
$\sigma$--algebra $\mathcal{B}(\M)$,  if the inverse image of any
set in $\F$ or in $\mathcal{B}(\M)$  belongs to $\F$. The latter
will be also called {\em random variable} with values in $\M$.
\par

We will be particulary interested in the case where the range of a
random variable is a {\em Polish  space}, namely  a complete and
separable metric space. By $\D C(\R^N)$ and
$\D{Lip}_\kappa(\R^n)$, we will denote the Polish
 space of continuous and
Lipschitz--continuous real functions (with Lipschitz constant less
than or equal to $\kappa>0$), defined in  $\R^N$, both endowed
with the metric $d$ inducing the topology of uniform convergence
on compact subsets of $\R^N$.
We will use the expressions {\em continuous random function}, {\em
$\kappa$--Lipschitz random function}, respectively, for the
previously introduced random variables. We will more simply say
{\em Lipschitz random function} to mean a  $\kappa$--Lipschitz
random function for some $\kappa>0$. See \cite{DavSic08}  for more
detail on this point.

We proceed by recalling some basic facts on convergence in
probability.  Given a Polish space $(\FF, d)$ and a sequence
$(f_n)_n$ of random variables taking values in $\FF$, we will say
that
 $f_n$ converge to $f$ {\em in probability} if, for every
$\eps>0$,
\[
\PP\left(\{\omega\in\Omega\,:\,d(f_n(\omega),f(\omega))>\eps\}\right)\to
0\quad\hbox{as $n\to +\infty.$}
\]
The limit $f$ is still a random variable. Since $\FF$ is a
separable metric space, almost sure convergence, i.e.
$d\left(f_n(\omega),f(\omega)\right)\rightarrow 0$ a.s. in $\omega$,
implies convergence in
probability, while the converse is not true in general. However,
the following characterization holds:

\begin{teorema}\label{tmea}
Let $f_n,f$ be random variables with values in $\FF$. Then $f_n\to
f$ in probability if and only if every subsequence $(f_{n_k})_k$
has  a  subsequence converging to $f$ a.s..
\end{teorema}

We denote  by $L^0(\Omega,\FF)$ the space made up by the
equivalence classes of random variables with value in $\FF$  for
the relation of  almost sure equality. For every $f,g\in
L^0(\Omega,\FF)$, we set
\[
\alpha(f,g):=\inf\{\eps\geq
0\,:\,\PP\big(\{\omega\in\Omega\,:\,d(f(\omega),g(\omega))>\eps\}\big)\leq\eps\}.
\]

\begin{teorema}\label{teo Ky Fan}
 $\alpha$ is a metric, named after Ky Fan, which metrizes
convergence in probability, i.e. $\alpha(f_n,f)\to 0$ if and only
if $f_n\to f$ in probability, and turns $L^0(\Omega,\FF)$ into a
Polish space.
\end{teorema}

An {\em $N$--dimensional dynamical system} $(\tau_x)_{x\in\R^N}$
is defined as a family of mappings $\tau_x:\Omega\to\Omega$ which
satisfy the following properties:
\begin{enumerate}
\item[{\em (1)}] the {\em group property:} $\tau_0=id$,\quad
$\tau_{x+y}=\tau_x\comp\tau_y$;

\item[{\em (2)}] the mappings $\tau_x:\Omega\to\Omega$ are
measurable and measure preserving, i.e. $\PP(\tau_x E)=\PP(E)$ for
every $E\in\F$;

\item[{\em (3)}] the map $(x,\omega)\mapsto \tau_x\omega$ from
$\R^N\times\Omega$ to $\Omega$ is jointly measurable, i.e.
measurable  with respect to the product $\sigma$--algebra
$\mathcal B (\R^N)\otimes\F$.
\end{enumerate}

We will moreover assume that $(\tau_x)_{x\in\R^N}$ is {\em
ergodic,} i.e. that one of the following equivalent conditions
hold:
\begin{itemize}
\item[{\em (i)}] every measurable function $f$ defined on $\Omega$
such that, for every $x\in\R^N$, $f(\tau_x\omega)=f(\omega)$ a.s.
in $\Omega$, is almost surely constant; \item[{\em (ii)}] every
set $A\in\F$ such that $\PP(\tau_x A\,\Delta\, A)=0$ for every
$x\in\R^N$ has probability either 0 or 1, where $\Delta$ stands for
the symmetric difference.
\end{itemize}

Given   a random variable $f:\Omega\to\R$, for any fixed
$\omega\in\Omega$ the function $x\mapsto f(\tau_x\omega)$ is said
to be a {\em realization of $f$.} The following properties follow
from Fubini's Theorem, see \cite{JiKoOl}:  if $f\in L^p(\Omega)$,
then $\PP$--almost all its realizations belong to
$L^p_{loc}(\R^N)$; if $f_n \rightarrow f$ in $L^p(\Omega)$, then
$\PP$--almost all realizations of $f_n$ converge to the
corresponding realization of $f$ in $L^p_{loc}(\R^N)$. The
Lebesgue spaces on $\R^N$ are understood with respect to the
Lebesgue measure.\par

The next lemma guarantees that a modification of a random variable
on a set of zero probability does not affect its realizations on
sets of positive Lebesgue measure on $\R^N$, almost surely in
$\omega$. The proof is based on Fubini's Theorem again, see Lemma
7.1 in \cite{JiKoOl}.

\begin{lemma}\label{lemmino utile}
Let $\widehat\Omega$ be a set of full measure in $\Omega$. Then
there exists a set of full measure
$\Omega'\subseteq\widehat\Omega$ such that for any
$\omega\in\Omega'$ we have $\tau_x\omega\in\widehat\Omega$ for
almost every $x\in\R^N$.
\end{lemma}

A jointly measurable function $v$ defined in $\R^N \times \Omega$
is said {\em stationary } if, for every $z \in \R^N$, there exists
a set $\Omega_z$ with probability $1$ such that for every
$\omega\in\Omega_z$
\[ v(\cdot + z, \omega)= v(\cdot, \tau_z \omega) \quad\text{on $\R^N$}
\]
It is clear that a real random
variable $\phi$  gives rise to a stationary function $v$
 by setting $v(x,\omega)=
\phi(\tau_x \omega)$. Conversely, according to Proposition 3.1 in
\cite{DavSic08}, a stationary function $v$ is, a.s. in $\omega$,
the realization of the measurable function $\omega \mapsto
v(0,\omega)$. More precisely,  there exists  a set $\Omega'$ of
probability $1$ such that for every $\omega\in\Omega'$
\begin{equation}\label{staz}
   v(x,\omega)=  v(0, \tau_x \omega) \quad\text{for a.e. $x\in\R^N$.}
\end{equation}

With the term {\em (graph--measurable ) random set} we indicate  a
set--valued function $X:\Omega\to\mathcal B(\R^N)$ with
\[
\Gamma(X):=\left\{(x,\omega)\in\R^N\times\Omega\,:\,x\in
X(\omega)\,\right\}
\]
jointly measurable in $\R^N\times\Omega$. A random set $X$ will be
qualified as {\em stationary} if for every  for every $z\in\R^N$,
there exists a set $\Omega_z$ of probability 1 such that
\begin{equation}\label{def stationary set}
X(\tau_z\omega)=X(\omega)-z\qquad\hbox{for every
$\omega\in\Omega_z$.}
\end{equation}

We use a stronger notion of measurability, which is usually  named
in the literature after Effros,  to define a {\em closed random
set}, say $X(\omega)$. Namely we require $X(\omega)$ to be a
closed subset of $\R^N$ for any $\omega$ and
\begin{equation*}\label{effros}
    \{\omega \, : \, X(\omega) \cap K \neq \emptyset \} \in \F
\end{equation*}
with $K$ varying  among the compact (equivalently, open) subsets
of $\R^N$. This condition can be analogously expressed by saying
that $X$ is measurable with respect to the Borel $\sigma$--algebra
related to the Fell topology on the family of closed subsets of
$\R^N$. This, in turn, coincides with the Effros
$\sigma$--algebra. If $X(\omega)$ is measurable in this sense then
it is also graph--measurable, see \cite{Molchanov} for more
details.

A closed random set $X$ is called stationary if it, in addition,
satisfies \eqref{def stationary set}. Note that in this event the
set $\{\omega\,:\,X(\omega)\not=\emptyset\,\}$, which is
measurable by the Effros measurability of $X$, is invariant with
respect to the group of translations $(\tau_x)_{x\in\R^N}$ by
stationarity, so it has probability either 0 or 1 by the
ergodicity assumption.

\begin{prop}\label{rando} Let $f$ be a continuous  random
function and $C$ a closed subset of $\R$. Then
\begin{equation*}\label{rando_1}
    X(\omega):= \{x \, :\, f(x,\omega)\in C\}
\end{equation*}
is a closed  random set in $\R^N$. If in addition $f$ is
stationary, then $X$ is stationary.
\end{prop}

See \cite{DavSic08} for a proof.

For a random stationary set $X$ it is immediate, by exploiting that
the maps $\{\tau_x\}_{x \in \R^N}$  are measure preserving, that
$\PP (X^{-1}(x))$ does not depend on $x$, where
\[X^{-1}(x)=\{\omega \, :\, x \in X(\omega)\}.\] Such quantity
 will be called {\em volume
fraction} of $X$ and denoted by $q_X$. Note that to any measurable
subset $\Omega'$ of $\Omega$ it can be associated a stationary set
$Y$  through the formula
\begin{equation*}
    Y(\omega):= \{x \,:\, \tau_x \omega \in \Omega'\}.
\end{equation*}
In this case $Y^{-1}(x)= \tau_{-x}\Omega'$, and so
$q_Y=\PP(\Omega')$.
By exploiting the ergodicity assumption and Birkhoff Ergodic
Theorem it is possible to derive an interesting  information on
the asymptotic structure of closed stationary sets.It says, in
particular,  they are spread with some uniformity in the space. We
refer the reader to \cite{DavSic08} for a proof.

\begin{prop}\label{prop stationary random set}
Let $X$ be an almost surely nonempty stationary closed random set in
$\R^N$. Then for every $\eps>0$ there exists $R_\eps>0$ such that
\[
\lim_{r\to +\infty} \frac{|\left(X(\omega)+B_R\right)\cap
B_r|}{|B_r|}\geq 1-\eps\qquad\hbox{a.s. in $\Omega$,}
\]
whenever $R\geq R_\eps$.
\end{prop}

\indent Given a Lipschitz random function $v$, we set
\[
\Delta_v(\omega):=\left\{\,x\in\R^N\,:\,v(\cdot,\omega)\ \hbox{is
differentiable at $x$}\, \right\}.
\]

\begin{definizione}\label{def v with st increments}
A random Lipschitz function $v$ is said to have
 {\em stationary increments} if, for every
    $z\in\R^N$, there exists a set $\Omega_z$ of probability 1
    such that
\begin{equation}\label{def stationary increments}
    v(x+z,\omega)-v(y+z,\omega)=v(x,\tau_z
    \omega)-v(y,\tau_z\omega)\quad\hbox{for all  $x,y\in\R^N$}
\end{equation}
    for every $\omega\in\Omega_z$.
\end{definizione}

    The following holds:

\begin{prop}\label{prop Delta_v random}
Let $v$ be a Lipschitz random function, then $\Delta_v$ is a
random set. In addition, it is stationary with volume
 fraction $1$ whenever $v$ has stationary increments.
\end{prop}

Let $v$ be a Lipschitz random function with stationary gradient.
For every fixed $x\in\R^N$, the random variable $Dv(x,\cdot)$ is
well defined on $\Delta_v^{-1}(x)$, which has probability 1 since
$\Delta_v$ is a stationary set with volume fraction 1.
Accordingly, we can define the mean $\EE(Dv(x,\cdot))$, which is
furthermore independent of $x$ by the stationary character of
$Dv$. In the sequel, we will be especially interested in the case
when this mean is zero.

\begin{definizione}\label{def ammissible}
A Lipschitz random function will be called {\em admissible} if it
has stationary increments and gradient with mean 0.
\end{definizione}

We state two characterizations of admissible random functions, and
a result that guarantees that stationary Lipschitz random
functions are admissible.

\begin{teorema}\label{teorema mean 0}
A Lipschitz random function $v$  with stationary increments has
gradient with vanishing mean if and only if it is almost surely
sublinear at infinity, namely
\begin{equation}\label{mean1}
    \lim_{|x|\to
    +\infty}\frac{v(x,\omega)}{|x|}=0\qquad\hbox{a.s. in $\omega$.}
\end{equation}
\end{teorema}
\smallskip
\begin{teorema}\label{teorema constant mean}
A Lipschitz random function  $v$ with stationary increments has
gradient  with vanishing mean if and only if
\begin{equation}\label{ciao}
   x\mapsto \EE(v(y,\cdot)- v(x,\cdot))=0 \quad\text{for any $x,y\in\R^N$.}\bigskip
\end{equation}
\end{teorema}

\begin{teorema}\label{teo admissible}
Any stationary Lipschitz random function $v$ is
admissible.\medskip
\end{teorema}

\indent Notice that the mean $\EE(v(x,\cdot))$ of a Lipschitz
random function  is independent of $x$, so when such a quantity is
finite Theorem \ref{teo admissible} is just a consequence of
Theorem \ref{teorema constant mean}.
\end{section}

\begin{section}{Stochastic Hamilton--Jacobi equations}\label{sez HJ}

We consider an Hamiltonian
\[
H:\R^N\times\R^N\times\Omega\to\R
\]
satisfying the following conditions:

\begin{itemize}
    \item[(H1)] the map $\omega\mapsto H(\cdot,\cdot,\omega)$
    from $\Omega$ to the Polish space $C(\R^N\times\R^N)$ is
    measurable;\smallskip
    \item[(H2)] for every $(x,\omega)\in\R^N\times\Omega$,
    $\ H(x,\cdot,\omega)\ \hbox{is  convex on $\R^N$;}$ \smallskip
     \item[(H3)]  there
     exist two superlinear functions $\alpha,\beta:\R_+\to\R$ such
     that
     \[
     \alpha\left(|p|\right)\leq H(x,p,\omega)\leq \beta\left(|p|\right)\qquad\hbox{for all
     $(x,p,\omega)\in\R^N\times\R^N\times\Omega$;}
     \]
      \item[(H4)] for every $(x,\omega)\in\R^N\times\Omega$,
      the set of minimizers of $H(x,\cdot,\omega)$
    has empty interior;\smallskip
    \item[(H5)] $H(\cdot+z,\cdot,\omega)=H(\cdot,\cdot,\tau_z\omega)$ for
    every $(z,\omega)\in\R^N\times\Omega$.
\end{itemize}

\begin{oss}
Condition (H3) is equivalent to saying that $H$ is superlinear and
locally bounded in $p$, uniformly with respect to $(x,\omega)$. We
deduce from (H2)
\begin{equation}\label{lippo}
  | H(x,p,\omega)-H(x,q,\omega)| \leq L_R |p-q| \quad\text{for all
  $x$, $\omega$, and $p$, $q$ in $B_R$},
\end{equation}
where $L_R=\sup\{\, |H(x,p,\omega)|\,:\,(x,\omega) \in \R^N \times
\Omega, \,|p|\leq R+2\,\},$ which is finite thanks to (H3). For a comment
on hypothesis (H4), see Remark \ref{oss h4}.
\end{oss}

\begin{oss}\label{oss almost-periodic}
Any given periodic, quasi--periodic or almost--periodic
Hamiltonian $H_0:\R^N\times\R^N\to\R$ can be seen as a specific
realization of a suitably defined stationary ergodic Hamiltonian,
cf. Remark 4.2 in \cite{DavSic08}. In the periodic and
quasi--periodic cases we take as $\Omega$ a $k$--dimensional
torus, with $k$ suitably chosen, which is separable both from the
topological and the measure theoretic viewpoint. In the
almost--periodic case, the usual construction is to take as
$\Omega$ the Bohr compactification of $\R^N$, which however is not
separable, cf. \cite{AmFrid}. In order to include this interesting
case in our treatment, we will show in the Appendix that, for a
given almost--periodic Hamiltonian $H_0$, it is possible to
construct a separable probability space $\Omega$, equipped with an
ergodic group of translations, such that $H_0$ can be seen as a
specific realization of a stationary ergodic Hamiltonian.
\end{oss}

For every $a\in\R$, we are interested  in the stochastic
Hamilton--Jacobi equation
\begin{equation}\label{eq HJa}
H(x,Dv(x,\omega),\omega)=a\qquad\hbox{in $\R^N$.}
\end{equation}
The material we are about to expose has been already presented in
\cite{DavSic08, DS1-08}, to which we refer for the details. Here
we just recall the main items.
\par

We say that a Lipschitz random function is a {\em solution} (resp.
{\em subsolution}) of \abra{eq HJa} if it is a viscosity solution
(resp.  a.e. subsolution) a.s. in $\omega$ (see \cite{BCD97, Ba94}
for the definition of viscosity (sub)solution in the deterministic
case). Notice that any such subsolution is almost surely in
$\D{Lip}_{\kappa_a}(\R^n)$, where
\begin{equation}\label{def kappa_a}
\kappa_a:=\sup\{\,|p|\,:\,H(x,p,\omega)\leq a\ \hbox{for some
$(x,\omega)\in\R^N\times\Omega$}\,\},
\end{equation}
which is finite thanks to (H3). We are interested in the class of
{\em admissible subsolutions}, hereafter denoted by $\S_a$, i.e.
random functions taking values in $\D{Lip}_{\kappa_a}(\R)$ with
stationary increments and zero mean gradient that are subsolutions
of \eqref{eq HJa}. An admissible solution will be also named {\em
exact corrector}, remembering its role in homogenization. Further,
for any $\delta >0$, a random function $v_\delta$ will be called a
{\em $\delta$--approximate corrector} for the equation \eqref{eq
HJa} if it belongs to $\S_{a+\delta}$ and satisfies the
inequalities
 \[
a - \delta \leq H(x, Dv_\delta(x,\omega),\omega) \leq a +\delta
\]
in the viscosity sense a.s. in $\omega$. We say that \abra{eq HJa}
has {\em approximate correctors} if it admits
$\delta$--approximate correctors for any $\delta >0$. \par

%

\medskip

We proceed by defining the {\em free} and the {\em stationary
critical value}, denoted by $c_f(\omega)$ and $c$ respectively, as
follows:
\begin{eqnarray}
 c_{f}(\omega)&=& \inf\left\{ a\in\R\,:\, \text{\abra{eq HJa} has a subsolution
$v\in\D{Lip}(\R^N)$} \right\}, \label{def cf}\\
c &=& \inf\{a\in\R\,:\,\S_a\not =\emptyset\,\}.  \label{def c}
\end{eqnarray}
We emphasize that in definition \abra{def cf} we are considering
{\em deterministic} a.e. subsolutions $v$ of the equation \abra{eq
HJa}, where $\omega$ is treated as a fixed parameter. Furthermore,
we note that $c_f(\tau_z\omega)=c_f (\omega)$ for every
$(z,\omega)\in\R^N\times\Omega$, so that, by ergodicity, the
random variable $c_f(\omega)$ is almost surely equal to a
constant, still denoted by  $c_f$. Hereafter we will write
$\Omega_f$ for the set of probability $1$ where $c_f(\omega)$
equals $c_f$.

Concerning the definition of the critical value $c$, we notice
that the set appearing at the right--hand side of \eqref{def c} is
non void, since it contains the value $\sup_{(x,\omega)}
H(x,0,\omega)$, which is finite thanks to (H3). Moreover, the
infimum is attained. In fact, see
\cite{DS1-08, LiSou03}

\begin{teorema}
$\S_c\not =\emptyset.$
\end{teorema}

It is apparent by the definitions that $c \geq c_f$. A more
precise result, establishing the relation with the effective
Hamiltonian obtained via the homogenization \cite{ReTa00,
Souga99}, will be discussed in the next section.


In the sequel, we mostly focus our attention on the {\em critical
equation}
\begin{equation}\label{eq critica}
H(x,Dv(x,\omega),\omega)=c\qquad\hbox{in $\R^N$.}
\end{equation}
The relevance of the critical value is given by the following
result, see Theorem 4.5 in \cite{DavSic08} for the proof.

\begin{teorema}\label{correctors}
The critical equation \abra{eq critica}  is the unique among the
equations \abra{eq HJa} for which either an exact corrector or
approximate correctors may exist.
\end{teorema}

\smallskip

Following the so called metric method for the analysis of \abra{eq
HJa}, see \cite{FaSic03}, we introduce an intrinsic
 path distance.  In next formulae
we assume that $a \geq c_f$ and  $\omega \in \Omega_f$. We  start
by defining  the sublevels
\[
Z_a(x,\omega):=\{p\,:\,H(x,p,\omega)\leq a\,\},
\]
and the related support functions
    \begin{equation*}\label{sigma}
    \sigma_a(x,q,\omega):=\sup\left\{\langle q,p\rangle\,:\,p \in Z_a(x,\omega)\,\right\}.
    \end{equation*}
It comes from \abra{lippo} (cf. Lemma 4.6 in \cite{DavSic08})
that, given $b>a$, we can find $\delta=\delta(b,a)>0$ with
\begin{equation}\label{spessore}
Z_a(x,\omega)+B_\delta\subseteq Z_b(x,\omega) \qquad\hbox{for
every $(x,\omega)\in\R^N\times\Omega_f$.}
\end{equation}
This property is needed in the proof of Theorem \ref{correctors}.
It is straightforward to check that $\sigma_a$ is convex in $q$,
upper semicontinuous in $x$ and, in addition, continuous whenever
$Z_a(x,\omega)$ has nonempty interior or reduces to a point. We
extend the definition of $\sigma_a$ to
$\R^N\times\R^N\times\Omega$ by setting
$\sigma_a(\cdot,\cdot,\omega)\equiv 0$ for every
$\omega\in\Omega\setminus\Omega_f$. With this choice, the function
$\sigma_a$ is jointly measurable in $\R^N \times\R^N\times \Omega$
and enjoys the stationarity property
\[
\sigma_a(\cdot+z,\cdot,\omega)=\sigma_a(\cdot,\cdot,\tau_z\omega)\quad\hbox{for
every $z\in\R^N$ and $\omega\in\Omega$.}
\]
We define the semidistance $S_a$  as
    \begin{equation}\label{eq S}
    S_a(x,y,\omega)=\inf\left\{\int_0^1 \sigma_a(\gamma(s),\dot\gamma(s),\omega)\,\dd
    s\,:\, \gamma\in\D{Lip}_{x,y}([0,1],\R^N)\, \right\},
    \end{equation}
The function $S_a$ is measurable on $\R^N\times\R^N\times\Omega$
with respect to the product $\sigma$--algebra $\mathcal
B(\R^N)\otimes\mathcal B(\R^N)\otimes\F$, and satisfies the
following properties:
\begin{eqnarray*}
S_a(x,y,\tau_z\omega)&=&S(x+z,y+z,\omega)\\
S_a(x,y,\omega)&\leq& S_a(x,z,\omega)+S_a(z,y,\omega)\\
S_a(x,y,\omega)&\leq& \kappa_a|x-y|
\end{eqnarray*}
for all $x,y,z\in\R^N$ and $\omega\in\Omega$. 


%
In the study of equation \eqref{eq HJa}, a special role is played
by the classical {\em (projected) Aubry set}  (cf.
\cite{FaSic03}), defined for every $\omega\in\Omega_f$ as the
collection of points $y\in\R^N$ such that
\begin{equation*}\label{eq A}
\inf\left\{\int_0^1\sigma_{a}(\gamma,\dot\gamma,\omega)\,\dd
s\,:\,\gamma\in\D{Lip}_{y,y}([0,1],\R^N),\,\hh^1(\gamma)\geq\delta
\right\}=0
\end{equation*}
{for some $\delta>0$, }or, equivalently (cf. \cite[Lemma
5.1]{FaSic03}), for any $\delta>0$.
From the Aubry--Mather theory for deterministic Hamiltonians we
know that, when $a>c_f$, this set is empty for all
$\omega\in\Omega_f$, i.e. almost surely. Hence, the only
interesting case is the one corresponding to $a=c_f$. Hereafter we
will denote by $\A_{f}(\omega)$ the collection of points $y$ of
$\R^N$ enjoying the above condition with $a=c_f$. The set
$\A_{f}(\omega)$ is closed for every $\omega\in\Omega$.

\par

We will also use later an equivalent definition of $\A_f(\omega)$,
see \cite{CoIt00}. For every $\omega\in\Omega$, let
\[
L(x,q,\omega)
    :=
   \max_{p\in\R^N}\left\{\langle p,q\rangle -
   H(x,p,\omega)\right\},\qquad\hbox{$(x,q)\in\R^N\times\R^N$}
\]
and, for every $t>0$,
\[
h_t(x,y,\omega):=\inf \left\{\int_0^t
\left(L(\gamma,\dot\gamma,\omega)+c\right)\,\dd
s\,:\,\gamma(0)=x,\,\gamma(t)=y\,\right\},
\qquad\hbox{$x,\,y\in\R^N$}.
\]
Then
\begin{equation}\label{aub}
    \A_f(\omega)=\{\,y\in\R^N\,:\ \liminf_{t\to +\infty}
    h_t(y,y,\omega)=0\,\}.
\end{equation}

\medskip

 In the next theorem we outline the main
deterministic properties linking $\A_{f}(\omega)$ to equation
(\ref{eq HJa}), see \cite{FaSic03}.

\begin{teorema}\label{teo 6.7}
Let $\omega\in\Omega_f$. The following holds:
\begin{itemize}
\item[\em (i)] Assume that $\A_{f}(\omega)\not=\emptyset$. If
$w_0$ is a function defined on $C \subset \A_{f}(\omega)$
    such that
    \[
    w_0(x)-w_0(y)\leq S_{c_f}(y,x,\omega)\quad\hbox{for every $x,y\in C$,}
    \]
    then the function
    \begin{equation*}
    w(x):=\min_{y\in C} \big(w_0(y)+S_{c_f}(y,x,\omega)\big)\qquad\hbox{$x\in \R^N$}
    \end{equation*}
    is the maximal subsolution of \abra{eq HJa} with $a=c_f$ equaling $w_0$ on $C$,
    and a solution as well.\medskip
\item[\em (ii)] Let $U$ be a bounded open subset of $\R^N$, and assume
that either $a>c_f$, or $a=c_f$ and
$U\cap\A_f(\omega)=\emptyset$. Let $w_0$ be a function defined on
$\partial U$ such that
\[
w_0(x)-w_0(y)\leq S_{a}(y,x,\omega)\quad\hbox{for every $x,y\in
\partial U$.}
\]
Then the function
\[
w(x):=\inf_{y\in \partial U}
\big(w_0(y)+S_{a}(y,x,\omega)\big)\qquad\hbox{$x\in U$}
\]
is the unique  viscosity solution of the Dirichlet Problem:
\begin{eqnarray*}
\begin{cases}
H(x,D \phi(x),\omega)=a&\qquad\hbox{in $U$}\\
\phi(x)=w_0(x)&\qquad\hbox{on $\partial U$.}
\end{cases}
\end{eqnarray*}
\item[\em (iii)] Assume that $a=c_f$ and  let $U$ be a bounded open subset of $\R^N$
with  $U\cap\A_f(\omega) \neq \emptyset$. Let $w_0$ be  a function
defined in $\partial U \cup \A_f$ $1$--Lipschitz continuous with
respect to $S_a$. Then the function
\[
w(x):=\inf \big \{w_0(y)+S_{a}(y,x,\omega)\,:\,y\in \partial U
\cup (U \cap \A_f)\big\}\qquad\hbox{$x\in U \setminus \A_f$}
\]
is the unique  viscosity solution of the Dirichlet Problem:
\begin{eqnarray*}
\begin{cases}
H(x,D \phi(x),\omega)=a&\qquad\hbox{in $U \setminus \A_f$}\\
\phi(x)=w_0(x)&\qquad\hbox{on $\partial U \cup (U \cap \A_f)$.}
\end{cases}
\end{eqnarray*}
\end{itemize}
\end{teorema}

\bigskip

\indent We define, for every $\omega\in\Omega$, the {\em set of
equilibria},  as follows:
\begin{equation*}
\Eq(\omega):=\{y\in\R\,:\,\min_p H(y,p,\omega)=c_f\,\}.
\end{equation*}
The set $\Eq(\omega)$ is a (possibly empty)  closed subset of
$\A_{f}(\omega)$ (cf. \cite[Lemma 5.2]{FaSic03}).
It is apparent that  ${c_f} \geq \sup_{x \in \R^N}
\min_{p \in \R^N} H(x,p,\omega)$ a.s. in $\omega$;  we point out
that $\Eq(\omega)$ is nonempty if and only if the previous formula
holds with an equality. In this case, $\Eq(\omega)$ is made up by
the points $y$ where the maximum is attained.

\begin{oss}\label{oss h4}
The inclusion $\Eq(\omega)\subseteq \A_f(\omega)$
depends on the fact that the $c_f$--sublevel
$\{p\,:\,H(y,p,\omega)\leq {c_f}\}$ is non--void and has empty
interior when $y \in \Eq(\omega)$. The latter is a consequence of
(H4), and this is actually the unique point where such condition
is used.
\end{oss}

We recall for later use a result from \cite{DS1-08}.

\begin{prop} \label{equi}
$\Eq(\omega)$ and $\A_f(\omega)$ are stationary closed random
sets.
\end{prop}

\end{section}

\begin{section}{Lax formula and closed random sets}\label{sez Lax}

\indent In this section we give a stochastic version of  Lax
formula and investigate when it  provides an exact corrector.
\smallskip\\
\indent Let $C(\omega)$ be an almost surely
nonempty stationary closed random set in $\R^N$. Take  a Lipschitz
random function $g$ and set, for $a\geq c_f$,
\begin{equation}\label{def u}
u(x,\omega):= \inf\{g(y,\omega)+S_a(y,x,\omega)\,:\,y\in
C(\omega)\,\}\quad \hbox{$x\in\R^N$,}
\end{equation}
where we agree that $u(\cdot,\omega)\equiv 0$ when either
$C(\omega)=\emptyset$ or the above infimum is equal to $-\infty$.
The following holds, see \cite{DavSic08,DS1-08}:

\begin{prop}\label{prop Lax in S} Let $a\geq c_f$ and $C(\omega)$,
$u$ as above.
\begin{itemize}
\item[{\em (i)}] Let $g$ be a stationary random function and
assume that the infimum  in \eqref{def u} is a.s. finite. Then $u$
is a stationary random variable belonging to $\S_a$ and satisfies
$u(\cdot,\omega)\leq g(\cdot,\omega)$ on $C(\omega)$ a.s. in
$\omega$. Moreover, $u$ is a viscosity solution of \eqref{eq HJa}
in $\R^N\setminus {C(\omega)}$ a.s. in $\omega$.\smallskip

\item[{\em (ii)}] Assume $g \in\S_a$. Then the random Lipschitz
function  $u$ belongs to $\S_a$ and satisfies
$u(\cdot,\omega)=g(\cdot,\omega)$ on $C(\omega)$ a.s. in $\omega$.
Moreover, $u$ is a viscosity solution of \eqref{eq HJa} in
$\R^N\setminus {C(\omega)}$ a.s. in $\omega$.
\end{itemize}
\end{prop}

We recall that the effective Hamiltonian $\overline H$ is the
function associating to any $P \in \R^N$ the critical value of the
Hamiltonian $H(x, p+P,\omega)$, equivalently it can be defined by
homogenization, see \cite{ReTa00, Souga99}. It can be proved, see
\cite{DavSic08,DS1-08}, that it is convex and superlinear, and
$\min_{\R^N} \overline H=c_f$. For any $a \geq c_f$   we denote by
$\overline Z_a$  the $a$--sublevel of $\overline H$.  By making
use of Propositions \ref{prop Lax in S} and \ref{equi}, the
following result  has been proved in \cite{DS1-08}

\begin{teorema}\label{nesiste}\ \vspace{1ex}
\begin{itemize}
    \item[{\em (i)}] If $c=c_f$ and the classical Aubry set $\A_f(\omega)$ is almost
    surely nonempty, then the extension of any $g \in \S_c$ from
    $\A_f$ through Lax formula with distance $S_c$ provides an exact corrector for
    \eqref{eq critica};
    \smallskip
    \item[{\em (ii)}] If  $0 \in \mathrm{Int}\big(\overline{Z}_c\big)$, then $c=c_f$ and there
    exists an exact  corrector for \eqref{eq critica} if and
    only if  the classical Aubry set $\A_f(\omega)$ is almost surely
    nonempty.\medskip
\end{itemize}
\end{teorema}

  To ease notations, from now on we will always write
$S$, $\sigma$ and $\S$ in place of $S_c$, $\sigma_c$ and $\S_c$,
respectively.\smallskip\par

The next result shows that the property that the Lax formula with
source a random set $C(\omega)$ and trace $g \in \S$ gives an
exact corrector can be solely detected looking at the behavior of
$g$ on $C$. This will be used in the next section for studying the
random Aubry set.

\begin{teorema}\label{teo cond C}
 Let $C(\omega)$ and $g$ be a stationary  closed random
set and a critical subsolution, respectively. Assume that either
$c > c_f$, or $c=c_f$ and $\A_{f}(\omega)\cap C(\omega)=
\emptyset$ a.s. in $\omega$.  Then the Lax extension of $g$ from
$C(\omega)$  with distance $S$, denoted by $u$, is an exact
corrector  if and only if
\begin{eqnarray}\label{cond C}
&\hbox{for any $y_0\in C(\omega)$ there exists a diverging
sequence
$(y_n)_n$  in $C(\omega)$}&\nonumber\\
&\hbox{such that}&\\
&\displaystyle g(y_0,\omega)=\lim_n\,
g(y_n,\omega)+S(y_n,y_0,\omega),&\nonumber
\end{eqnarray}
a.s. in $\omega$.
\smallskip
\end{teorema}

\begin{dimo}
\quad{\em \eqref{cond C} holds $\Rightarrow$ $u$ is an exact
corrector}\smallskip\\
\indent  In view of Proposition \ref{prop Lax in S}, we can select
a subset $\Omega'$  of $\Omega$ with $\PP(\Omega')=1$  such that
$C(\omega)\not =\emptyset$, \eqref{cond C} holds and
$u(\cdot,\omega)$ is a viscosity solution to \eqref{eq critica} in
$\R^N \setminus C(\omega)$, whenever $\omega \in \Omega'$. Let us
fix $\omega$ in $\Omega'$. If $u(\cdot,\omega)$ is not a critical
solution, there exist $x_0\in C(\omega)$ and a strict $C^1$
subtangent $\varphi$ to $u(\cdot,\omega)$ at $x_0$ with
\[
H(x_0,D\varphi(x_0),\omega)<c.
\]
By the usual technique  of pushing up such  test function,  we can
construct a  deterministic subsolution $v$  to
$H(x,Du,x,\omega)=c$ such that
\[
v(x_0)>u(x_0,\omega) \quad\hbox{and}\quad v(y_n)=u(y_n,\omega)
\quad\hbox{definitively in $n$.}
\]
For  $n$ sufficiently large  we  then get
\[
v(y_{ n})+S(y_{ n},x_0,\omega)<v(x_0),
\]
which is impossible by the subsolution property of  $v$.\medskip\\
\indent{\em $u$ is an admissible
solution $\Rightarrow$ \eqref{cond C} holds  }\smallskip\\
Let us fix $\omega\in\Omega$ such that $C(\omega)\not=\emptyset$,
$C(\omega) \cap \A_f(\omega) \neq \emptyset$, and
$u(\cdot,\omega)$ and $g(\cdot,\omega)$ are an admissible critical
solution and subsolution, respectively. These properties hold in a
subset of $\Omega$ with probability $1$. We introduce a partial
order relation in $C(\omega)$ by setting
\[
y_1\succ y_2 \quad \iff \quad
g(y_2,\omega)=g(y_1,\omega)+S(y_1,y_2,\omega).
\]
 We exploit the triangle inequality on $S$ and  the fact that
$g(\cdot,\omega)$ is a subsolution,  to see that this relation
enjoys the transitivity property. To prove that it is also
antisymmetric, we consider $y_1$, $y_2$ with $  y_1\succ y_2$ and
$ y_2\succ y_1$,
 accordingly
\begin{eqnarray*}
g(y_2,\omega)=g(y_1,\omega)+S(y_1,y_2,\omega)\quad\hbox{and}\quad
g(y_1,\omega)=g(y_2,\omega)+S(y_2,y_1,\omega).
\end{eqnarray*}
By summing up,  we get $S(y_1,y_2,\omega)+S(y_2,y_1,\omega)=0$,
which gives $y_1=y_2$, as desired,  since $\A_f(\omega) \cap
C(\omega)= \emptyset$.

For a fixed $y_0\in C(\omega)$, let
\[
C_{y_0}(\omega)=\{y\in C(\omega)\,:\,y\succ y_0\,\}.
\]
Using the continuity of $g(\cdot,\omega)$, $S(\cdot,y_0,\omega)$
and the closed character of $C(\omega)$, it is easy to check that
this set is closed.  If we show that $C_{y_0}(\omega)$ is
unbounded, the assertion is obtained.

Let us then assume, for purposes of contradiction, that
$C_{y_0}(\omega)$ is compact. We show that in this case
$C_{y_0}(\omega)$ admits a maximal element with respect to
$\succ$. Thanks to Zorn lemma, it suffices to prove:\smallskip\\
{\em Claim: any totally ordered subset $E$ of $C_{y_0}(\omega)$
admits an upper bound in $C_{y_0}(\omega)$.\smallskip\\}
\indent We first show that $\overline E$ is totally ordered, i.e.
$y\succ y'$, $y'\succ y$ or $y=y'$ for any pair $y$, $y'$ of
elements of $\overline E$. Let
\[
y=\lim_n y_n,\ \  y'=\lim_n y_n',\quad\hbox{with $y_n,y'_n\in E$
for every $n\in\N$.}
\]
If  definitively $y_n\succ y_n'$, then passing to  the limit in
the equality
\[
g(y_n',\omega)=g(y_n,\omega)+S(y_n,y_n',\omega)
\]
we get  $y\succ y'$. Similarly  $y'\succ y$ if $y_n'\succ y_n$
definitively. Finally,  if there exist two subsequences with
\[
y_{n_j}'\succ y_{n_j}\ \hbox{and}\ y_{n_k}\succ y_{n_k}',
\]
 we get  $y=y'$, for both $y\succ y'$ and $y'\succ
y$ hold, and $\succ$ enjoys the antisymmetric property.

Since  $\overline E$ is  compact, for every $\eps>0$ we find a
finite number $m=m(\eps)$ of points $y^\eps_1,\dots,y^\eps_m$ in
$\overline E$ such that $\overline E\subset \cup_i
B_\eps(y^\eps_i)$. Up to a reordering, we can as well assume
 $y^\eps_1 \succ y^\eps_j$ for all $j\not=1$. For every $y\in
\overline E$ and a suitable  $i\in\{1,\dots,m\}$ we have
\begin{eqnarray*}
    g(y,\omega)\!\!\!\!&\geq&\!\!\!\! g(y^\eps_i,\omega)+S(y^\eps_i,y,\omega)-2\kappa_c\,\eps\\
    \!\!\!\!&\geq&\!\!\!\!
    g(y^\eps_1,\omega)+S(y^\eps_1,y^\eps_i,\omega)+S(y^\eps_i,y,\omega)-2\kappa_c\eps
    \geq
    g(y^\eps_1,\omega)+S(y^\eps_1,y,\omega)-2\kappa_c\eps.\phantom{andrea}
\end{eqnarray*}
Taking  the limit as $\eps\to 0$  of $y^\eps_1$ and using the
compactness of $\overline E$,  we get an upper bound  for $E$, as
it was claimed.

We denote by  $\widetilde y$ a maximal element in
$C_{y_0}(\omega)$ with respect to $\succ$.   Since
$u(\cdot,\omega)$ agrees with $g(\cdot,\omega)$ on $C(\omega)$ and
is a viscosity solution of \eqref{eq critica}, and $\A_f(\omega)
\cap C(\omega) = \emptyset$, Theorem \ref{teo 6.7} yields that
there exists $y' \neq \widetilde y$ with
\[
g(\widetilde y,\omega)=u(y',\omega)+S(y',\widetilde y,\omega).
\]
If
\[
u(y',\omega)=g(z,\omega)+S(z,y',\omega) \qquad\hbox{for some $z\in
C(\omega)$,}
\]
then
\[
g(\widetilde y,\omega)=g(z,\omega)+S(z,y',\omega)+S(y',\widetilde
y,\omega) \geq g(z,\omega)+S(z,\widetilde y,\omega),
\] which, in turn, implies $z=\widetilde y$ since $z \in C_{y_0}(\omega)$ by the transitivity property of
$\succ$ and  $\widetilde y$ is maximal;   consequently
$S(\widetilde y,y',\omega)+S(y',\widetilde y,\omega) =0$ which is
impossible since $\widetilde y \neq y'$ and $\A_f(\omega) \cap
C(\omega) = \emptyset$. Therefore
\[
u(y',\omega)=\lim_n \big(\, g(y_n,\omega)+S(y_n,y',\omega)\big),
\]
for some diverging sequence $(y_n)_n$ in $C(\omega)$.  We derive
\begin{eqnarray*}
g(\widetilde y,\omega)&=&\lim_n\,\big(
g(y_n,\omega)+S(y_n,y',\omega)+S(y',\widetilde
y,\omega)\big)\\
&\geq& \lim_n \big(g(y_n,\omega)+S(y_n,\widetilde
y,\omega)\big)\geq g(\widetilde y,\omega),
\end{eqnarray*}
and, since $\widetilde y\succ y_0$,
\begin{eqnarray*}
  g(y_0,\omega)&=&g(\widetilde y,\omega)+S(\widetilde y,y_0,\omega) =
\lim_n g(y_n,\omega)+S(y_n,\widetilde y,\omega) + S(\widetilde
y,y_0,\omega) \\
  &\geq& \lim_n g(y_n,\omega)+S(y_n,y_0,\omega).
\end{eqnarray*}
Since the converse inequality also holds as $g(\cdot,\omega)$ is a
critical subsolution, we finally obtain that $y_n \in
C_{y_0}(\omega)$ for any $n$, which is impossible since $y_n$ is a
diverging sequence and $C_{y_0}(\omega)$  is a compact set, by
assumption.
\end{dimo}

\medskip

We point out that, in the previous theorem,  the argument for
deriving from \eqref{cond C}   that $u$ is an exact corrector, can
be used to get  a slight more general assertion, that we write
down below for later use.

\begin{cor}\label{cor cond C} Let $C(\omega)$ and $u$ be a stationary  closed random
set and the Lax extension of some critical subsolution from
$C(\omega)$ with distance $S$, respectively.  If for any $y_0\in
C(\omega)$ there exists  $y_1 \neq y_0$  with
\[ u(y_0,\omega)= u(y_1,\omega)+S(y_1,y_0,\omega),\]
a.s. in $\omega$, then  $u$ is an exact corrector.
\end{cor}

\medskip

We derive a further corollary of Theorem \ref{teo cond C}:

\begin{cor}\label{cor2 cond C}
Let $C(\omega)$ and $g$ be a stationary  closed random and  an
admissible critical subsolution, respectively. Assume that either
$c > c_f$, or $c=c_f$ and $\A_{f}(\omega)\cap C(\omega)=
\emptyset$ a.s. in $\omega$.  If the Lax extension of $g$ from
$C(\omega)$ with distance $S$  is an exact corrector  then
\begin{eqnarray}\label{cond2 C}
&\hbox{for any $y_0\in C(\omega)$ there exists a diverging
sequence
$(z_n)_n$  in $\R^n$}&\nonumber\\
&\hbox{such that}&\\
&\displaystyle g(y_0,\omega)= g(z_n,\omega)+S(z_n,y_0,\omega)
\quad\text{for any $n$},&\nonumber
\end{eqnarray}
a.s. in $\omega$.
\smallskip
\end{cor}

\begin{dimo}  Given $\omega$ in a subset of $\Omega$ with probability $1$ and $y_0 \in C(\omega)$,
there is, by Theorem \ref{teo cond C},  a diverging sequence
$(y_n)_n$ in  $C(\omega)$ satisfying \eqref{cond C}. Given $k \in
\N$, we can assume, without loss of generality, that $|y_n| > k$,
for any $n$. Let $(\xi_n)_n$ a  sequence of curves, defined in
$[0,1]$, joining $y_n$ to $y_0$ with
\begin{equation}\label{C1}
    \int_0^1\sigma(\xi_n,\dot\xi_n,\omega)\,\dd
s + g(y_n,\omega) \leq g(y_0,\omega) + 1/n \quad\text{for any $n
\in \N$.}
\end{equation}
Since $|y_n| > k$, there is, for any $n$, $t_n \in [0,1]$ with
\[|\xi(t_n)| = k.
\]
From \eqref{C1} we derive
\begin{eqnarray*}
S(y_n,\xi_n(t_n),\omega)&+&S(\xi_n(t_n),y_0,\omega)\\
    &\leq&
    \big( g(\xi_n(t_n),\omega)-g(y_n,\omega)\big)
    +
    \big( g(y_0,\omega)-g(\xi_n(t_n),\omega)\big)
    +
    1/n\qquad
\end{eqnarray*}
and taking into account  that $g$ is a critical subsolution, we
get
\[
\lim_n S(\xi_n(t_n), y_0, \omega)=\lim_n \
g\left(y_0,\omega\right) -g\left(\xi_n(t_n),\omega\right)
\]
For any limit point $z_k$ of
$(\xi_n(t_n))$, we find
\[g(y_0,\omega) -g(z_k,\omega) = S(z_k,y_0,\omega) \quad\text{where}\quad
|z_k|=k,
\]
and since $k \in \N$ was arbitrarily chosen, the assertion
follows.
\end{dimo}

\bigskip

\end{section}

\begin{section}{Random Aubry set}\label{sez Aubry}

We start by introducing a notion of Aubry set adapted to the
stationary ergodic setting, see also Remark 6.9 in
\cite{DavSic08}. To motivate it, we  recall that in the
deterministic case the Aubry set can be characterized by the
property that the critical intrinsic distance from any of its
points is a critical solution. Roughly speaking, the idea
underlying the next definition is to replace points by random
stationary closed subsets and make use of the Lax formula taking
as trace any critical admissible subsolutions.

\begin{definition}\label{def Aubry} A stationary closed random set $\A(\omega)$ is called
{\em random Aubry set} if
\begin{itemize}
    \item[\em (i)] the extension of any  admissible critical subsolution
from  $\A(\omega)$ via the Lax formula \eqref{def u} yields an
exact corrector;
    \item[\em (ii)] any  closed  random stationary set $C(\omega)$ enjoying the
previous property is almost surely contained in $\A(\omega)$.
\end{itemize}
\end{definition}

\medskip

We also need the following

\begin{definition} An admissible critical subsolution is called
{\em weakly strict} on some random set $X(\omega)$ if a.s. in
$\omega$
\[
v(x,\omega)-v(y,\omega)< S(y,x,\omega)\qquad\hbox{for every
$x,y\in X(\omega)$ with $x\not= y$.}
\]
\end{definition}

\medskip

The main result of the first part of the  section is

\begin{teorema}\label{teo stretta} Assume that  $c > c_f$ or $c=c_f$ and
$\A_f(\omega)= \emptyset$ a.s. in $\omega$. Then there exists a
critical admissible subsolution which is weakly strict in $\R^N
\setminus \A(\omega)$ a.s. in $\omega$. \end{teorema}

This, in particular, implies the existence of a critical
admissible subsolution,  weakly strict on the whole $\R^N$, if the
random Aubry set is almost surely empty.

\smallskip
We postpone the proof after some preliminary analysis. When
$c=c_f$ it is clear by Theorem \ref{nesiste} that $\A_f(\omega)
\subseteq \A(\omega)$,  and this inclusion
 can be strict a.s. in $\omega$. This occurs even in the periodic setting. Albert
Fathi provided us with an example of a periodic Hamiltonian for
which $\A_f$ is empty, while, of course, $\A$ is not. In this
example, however, $0\in\partial\overline Z_{c_f}$. Actually we
have:

\begin{prop}\label{coincide} Assume that $0 \in
\mathrm{Int}\left(\overline Z_{c_f}\right)$ and, consequently,
that $c=c_f$. Then $\A(\omega)= \A_f(\omega)$ a.s. in $\omega$.
Moreover $\A_f(\omega)$ is a uniqueness set for \eqref{eq
critica}.
\end{prop}

\begin{dimo} If $\A_f=\emptyset$ a.s. in $\omega$, then $\A$ is
also almost surely empty since no correctors can exist by Theorem
\ref{nesiste} (ii). Let us assume that $\A_f(\omega) \neq
\emptyset$ a.s. in $\omega$, and,  in addition, for purposes of
contradiction, that $\A_f(\omega) \subsetneq\A(\omega)$ a.s in
$\omega$. We claim that, in this case, there exists a closed
random stationary  a.s. nonempty   set $C(\omega)$ with
\[ C(\omega) \subset \A(\omega) \quad\text{and} \quad C(\omega) \cap
\A_f(\omega) = \emptyset \quad\text{a.s. in $\omega$.}\] For this,
we denote by $f(x,\omega)$ the Euclidean distance of $x$ from
$\A_f(\omega)$, for any $x$, $\omega$ (with the convention that it
is equal to $- \infty$ whenever $\A_f(\omega)$ is empty) and, for
any $n \in \N$, and consider the random set
\[ C_n(\omega):= \A(\omega) \cap \{x\,:\, f(x,\omega) \geq 1/n\}.\]
We see that it is  closed stationary taking into account that $f$
is a stationary continuous random function, Proposition
\ref{rando}, and the fact that the intersection of two closed
random stationary sets inherits the same  property. If
$C_n(\omega)= \emptyset$ a.s. in $\omega$, for any $n$, then
\[\A (\omega) \subset \bigcap_n \,\{x\,:\,f(x,\omega) < 1/n\} =
\A_f(\omega) \quad\text{a.s. in $\omega$,}\] which is in contrast
with our assumption. Accordingly,  there exists $n_0$ with
$C_{n_0} \neq \emptyset$ a.s in $\omega$. The claim is proved by
taking $C=C_{n_0}$.\par

Let now $u$ be any critical admissible subsolution. By the very
definition of random Aubry set, the Lax extension of $u$ from $C$
via $S$ yields an exact corrector, then, according to Theorem
\ref{teo cond C} and \eqref{cond C},   we find a.s. in $\omega$
\[u(y_0,\omega)=\lim_n u(y_n,\omega)+ S(y_n,y_0,\omega)\]
for any $y_0 \in C(\omega)$ and some diverging sequence $(y_n)_n$
of elements of $C(\omega)$. On the other side, since $0 \in
\mathrm{int}\overline Z_c$, we have a.s. in $\omega$
\begin{equation}\label{coincide1}
\lim_{|y| \rightarrow + \infty}  u(y,\omega)+ S(y,y_0,\omega)=+
\infty \quad\text{for any $y_0 \in \R^N$,}
\end{equation} which yields a
contradiction.\par

Let us finally prove the asserted uniqueness property of $\A_f$.
Let $v$ be  an exact corrector, we  fix $\omega$ such that
$\A_f(\omega) \neq \emptyset$, $v(\cdot, \omega)$ is a solution to
$H(x,Du,\omega)=c$ and \eqref{coincide1} holds true. We consider
the sequence of Dirichlet problems.
\begin{eqnarray*}
\begin{cases}
 H(x,D u,\omega)=c & \hbox{in
$B_n \setminus \A_f(\omega)$}\\
u(x)=v(x,\omega) & \hbox{in $\partial B_n \cup (B_n \cap
\A_f(\omega)$}.
\end{cases}
\end{eqnarray*}
 According to Theorem \ref{teo 6.7}, we find, for any $n$, the
relation
\[
v(0,\omega):=\inf \big \{v(y,\omega)+S(y,0,\omega)\,:\,y\in
\partial B_n \cup (B_n \cap \A_f)\big\}.
\]
Letting $n$  go to infinity and taking into account
\eqref{coincide1}, we deduce the existence of $y_0 \in
\A_f(\omega)$  satisfying
\[ v(0,\omega)= v(y_0,\omega) + S(y_0,0,\omega).\]
By applying the previous argument to any $x \in \R^N$ in place of
$0$, we finally get
\[v(x,\omega)= \inf\{ v(y,\omega) + S(y,x,\omega)\,:\, x \in
\A_f(\omega),\}\] which says that any exact corrector is the Lax
extension of its trace on $\A_f(\omega)$ a.s. in $\omega$. This
ends the proof.
\end{dimo}

\bigskip

\indent We assume from now on that $c > c_f$ or $c=c_f$ and
$\A_f(\omega)= \emptyset$ a.s. in $\omega$. For any $v\in\S$, we
define
\begin{equation}\label{def r_v}
r_v(x,\omega)=\max\{r\geq 0\,:\,\inf_{y\in\partial
B_r(x)}\big(v(y,\omega)+S(y,x,\omega)\big)=v(x,\omega)\,
\},\quad\hbox{$\omega\in\Omega$.}
\end{equation}

\begin{prop}\label{prop r_v}
Let $v\in\S$. The following properties hold:
\begin{itemize}
    \item[{\em (i)}] the map $r_v:\R^N\times\Omega\to\R$ is jointly measurable
    in $\Omega \times \R^N$;\smallskip
    \item[{\em (ii)}] $r_v$ is stationary;\smallskip
    \item[{\em (iii)}] $r_v(\cdot,\omega)$ is upper
    semicontinuous on $\R^N$ for every $\omega\in\Omega$;\smallskip
    \item[{\em (iv)}] {for every $z\in\R^N$,}\quad
    $r_v(\cdot,\tau_z\omega)=r_v(\cdot+z,\omega)$\quad{a.s. in $\omega$;}\smallskip
    \item[{\em (v)}] if $v(\cdot,\omega)$ is the  local uniform
    limit in $\R^N$ of a sequence $v_n(\cdot,\omega)$, with
    $v_n\in\S$  for every $n$, then
    \[
    \limsup_{n\to +\infty}\, r_{v_n}(x,\omega)\leq
    r_v(x,\omega)\qquad\hbox{for every $x\in\R^N$.}
    \]
    \item[{\em (vi)}] if $\hat v:=v-v(0,\omega)$ for every $\omega$, then $r_{\hat v}= r_v$ in $\R^N\times\Omega$.
\end{itemize}
\end{prop}

\begin{dimo}
Let us denote by $\psi_r(x,\omega)$ the infimum appearing in
formula \eqref{def r_v}. Fix $r>0$ and let $(z_n)_n$ be a dense
subset of $\partial B_r$. It is clear that
\[
\psi_r(x,\omega)=\inf_{n\in\N}\left(v(x+z_n,\omega)+S(x+z_n,x,\omega)\right),
\]
which implies that $\psi_r$ is measurable on $\R^N\times\Omega$.
Let now $(r_n)_n$ be a dense subset of $\R_+$ and, for each
$n\in\N$, set
\[
E_n:=\{(x,\omega)\in\R^N\times\Omega\,:\,\psi_{r_n}(x,\omega)=v(x,\omega)\,\}.
\]
Then $r_v(x,\omega)=\sup_n r_n\,\cchi_{E_n}(x,\omega)$ on
$\R^N\times\Omega$, and this proves {\em(i)}. Assertions {\em
(ii)}--{\em (vi)}  follow from the very definition
of $r_v$ and the fact that $v$ has stationary increments.
\end{dimo}
\medskip\\

\medskip
We will also need the following:

\begin{lemma}\label{lemma utile}
Let $v\in\S$ and $\alpha>0$. Then the sets
\[
C_\alpha(\omega):=\{x\in\R^N\,:\,r_v(x,\omega)\geq
\alpha\,\},\quad
C_\infty(\omega):=\{x\in\R^N\,:\,r_v(x,\omega)=+\infty\,\}
\]
are stationary closed random sets.
\end{lemma}

\begin{dimo}
It is clear by Proposition \ref{prop r_v} {\em (ii)} that
$C_\alpha$ is stationary. In order to prove that $C_\alpha$ is a
closed random set, we note that
$C_\alpha(\omega)=\{x\in\R^N\,:\,G_\alpha(x,\omega)=v(x,\omega)\,\},$
where
\[
G_\alpha(x,\omega):=\min_{y\in\partial
B_\alpha}\big(v(x+y,\omega)+S(x+y,x,\omega)\,\big),\quad\hbox{$(x,\omega)\in\R^N\times\Omega$}.
\]
It is easily seen that $G_\alpha$ is jointly measurable and
continuous in $x$ for any fixed $\omega$, thus proving the
asserted property for $C_\alpha(\omega)$ in view of Proposition
\ref{rando}. The remainder of the statement follows since
$C_\infty(\omega)=\bigcap_n C_n(\omega)$ and the intersection of a
countable family of stationary closed random sets is still a
stationary closed random set.
\end{dimo}

\bigskip

Theorem \ref{teo cond C} suggests that the following identity should hold
\begin{equation*}
\A(\omega)  = \bigcap_{v\in\S} \{x\in\R^N \,:\, r_v(x,\omega) = +
 \infty\},
\end{equation*}
a.s in $\omega$. However this need not be true, the main
difficulty being  that the above intersection is not countable in
general. To avoid this problem, we essentially exploit the
separability assumption on $\Omega$. According to Theorems
\ref{tmea} and \ref{teo Ky Fan}, the family of renormalized
critical subsolution
\[
\widehat\S:=\{\hat v\in\S\,:\,\hat v(0,\omega)=0\ \hbox{for every
$\omega$}\}
\]
is a subspace of $L^0(\Omega,\D C(\R^N))$, in particular it is
separable with respect to the Ky Fan metric. Therefore there
exists a sequence of Lipschitz random functions  $(v_n)_n$ which
is dense in $\widehat S$ with respect to the convergence in
probability. That implies, in view of Theorem \ref{tmea}, that
$(v_n)_n$ is also dense for the almost sure convergence in $\D
C(\R^N)$. We have

\begin{teorema}\label{teo Aubry}
Let $(v_n)_n$ as above. Then
\begin{equation}\label{eq Aubry}
\A(\omega)  = \bigcap_{n\in\N} \{x\in\R^N \,:\, r_{v_n}(x,\omega)
= +
 \infty\}\qquad\hbox{a.s. in $\omega$.}
\end{equation}
\end{teorema}

\begin{dimo}
Let us denote by $C(\omega)$ the set appearing at the right--hand
side of \eqref{eq Aubry}, for every $\omega\in\Omega$. The fact
that $C(\omega)$ is a stationary closed random set follows from
Lemma \ref{lemma utile}. Let us show that $C(\omega)\subseteq
\A(\omega)$ a.s. in $\omega$. By definition of Aubry set, we need
to show that $C(\omega)$ enjoys item {\em (i)} in Definition
\ref{def Aubry}. According to Corollary \ref{cor cond C}, this
amounts to requiring   the following identity to hold almost
surely:
\begin{equation}\label{eq equivalente}
r_v(\cdot,\omega)\ > 0\quad\hbox{on $C(\omega)$,}
\end{equation}
whenever $v$ is the Lax extension of some admissible trace from
$C(\omega)$.  We set $\hat v=v-v(0,\omega)$ for every $\omega$.
Clearly $\hat v\in\widehat S$, so there exists a sequence
$(v_{n_k})_k$ and a set $\Omega_0$ of probability 1 such that
$v_{n_k}(\cdot,\omega)\ucv \hat v(\cdot,\omega)$ in $\R^N$ for
every $\omega\in\Omega_0$. By Proposition \ref{prop r_v} for any
such $\omega$ we get
\[
\limsup_{k\to +\infty} r_{v_{n_k}}(x,\omega)\leq r_{\hat
v}(x,\omega)=r_{v}(x,\omega)\qquad\hbox{for every $x\in\R^N$},
\]
thus proving \eqref{eq equivalente} by the definition of
$C(\omega)$. Conversely, since $\A(\omega)$ enjoys item {\em (i)}
in Definition \ref{def Aubry}, we have in particular that
$r_{v_n}(\cdot,\omega)\equiv +\infty$ on $\A(\omega)$ a.s. in
$\omega$ for every $n\in\N$  by Corollary \ref{cor2 cond C}. That
implies $\A(\omega)\subseteq C(\omega)$ and concludes the proof.
\end{dimo}

\bigskip

We proceed by showing the existence of a random function
$\overline v$ in $\S$ enjoying a minimality property.

\begin{prop} \label{stramin}
There exist  $\overline v \in \S$  such that, for every $v\in\S$,
the following inequality holds almost surely:
\[
 r_{\overline v}(x,\omega)\leq r_v(x,\omega)
 \quad\text{in $\R^N$.}
\]
In particular, $\A(\omega)=\{x\in\R^N \,:\, r_{\overline
v}(x,\omega) = + \infty\}$ a.s. in $\omega$.
\end{prop}

\begin{dimo}
Let us take a sequence of positive real numbers
$(\lambda_n)_n$ with $\sum_n \lambda_n=1$ and set
\begin{equation}\label{def overline v} \overline
v(x,\omega)=\sum_{n=1}^{+\infty} \lambda_n
v_n(x,\omega),\qquad\hbox{for every
$(x,\omega)\in\R^N\times\Omega,$}
\end{equation}
where $v_n$ are the renormalized critical subsolutions appearing in
\eqref{eq Aubry}. It is easy to check that $\overline v\in\S$. Let
$\Omega_0$ be a set of probability 1 such that for
$\omega\in\Omega_0$ all the functions $v_n(\cdot,\omega)$ are
subsolutions of the critical equation \eqref{eq critica}.  Let us
fix $\omega \in \Omega_0$ and $x \in \R^N$. If
$|y-x|>r_{v_n}(x,\omega)$ for some $n \in \N$, then
\begin{eqnarray*}
\overline v(x,\omega)-\overline v(y,\omega)&=&\sum_{k\not= n}
\lambda_k\left(v_k(x,\omega)-v_k(y,\omega)\right) +
\lambda_{n}\left(v_{n}(x,\omega)-v_{n}(y,\omega)\right)\\
&<& \sum_{k\not=n}\lambda_k S(y,x,\omega)+\lambda_{n}
S(y,x,\omega)=S(y,x,\omega).
\end{eqnarray*}
We derive
\begin{equation}\label{disug1}
    r_{\overline v}(x,\omega)\leq r_{v_n}(x,\omega) \quad\text{ for every
$x\in\R^N$ and $n\in\N$.}
\end{equation}
To show that \eqref{disug1} holds true almost surely when $v_n$ is
replaced by any $v\in\S$,  set $\hat v=v-v(0,\omega)$. Clearly
$\hat v\in\widehat S$, and being $(v_n)_n$ dense in $\widehat\S$
with respect to the almost sure convergence, we derive by
Proposition \ref{prop r_v} {\em (iii)}, {\em (vi)} that
\begin{equation*}
    \liminf_{n\to +\infty} r_{v_n}(x,\omega)
    \leq
    r_{\hat v}(x,\omega)
    =
    r_v(x,\omega)
    \qquad\hbox{for every $x\in\R^N$}
\end{equation*}
a.s. in $\omega$. This shows the minimality of $r_{\overline v}$ and, consequently,
that $\{x\in\R^N \,:\, r_{\overline v}(x,\omega)
= + \infty\}\subseteq \A(\omega)$ a.s. in $\omega$. The opposite inclusion holds as well
since $r_{\overline v}(\cdot,\omega)\equiv +\infty$
on $\A(\omega)$ a.s. in $\omega$ by definition of Aubry set,
as already remarked in the
proof of Theorem \ref{teo Aubry}.
\end{dimo}

\bigskip

{\bf\noindent Proof of Theorem \ref{teo stretta}.}\quad More
precisely, we will prove that  there exists $\overline v\in\S$
such that
\[
r_{\overline v}(\cdot,\omega)\equiv +\infty\quad\hbox{on $\A(\omega)$},
\qquad
r_{\overline v}(\cdot,\omega)\equiv 0\quad\hbox{in $\R^N\setminus \A(\omega)$},
\]
a.s. in $\omega$.\par

 Let $\overline v$ the admissible critical
subsolution given by Proposition \ref{stramin}, see \eqref{def
overline v}. If $\overline v$ is weakly strict on the whole
$\R^N$, we derive from Proposition \ref{stramin} that $\A(\omega)$
is almost surely empty and the assertion follows. If, on the other
hand, $\overline v$ is not weakly strict on $\R^N$, for a suitable
$\alpha>0$ the set
\[
C_\alpha(\omega):=\{ x \in \R^N \,:\, r_{\overline v}(x,\omega) \geq \alpha\}
\]
is a.s. nonempty, and is in addition a closed stationary random
set by Lemma \ref{lemma utile}. In view of Proposition
\ref{stramin}
\begin{equation}\label{eq inclusione}
C_\alpha(\omega)\supseteq\{x\in\R^N \,:\, r_{\overline v}(x,\omega)
= + \infty\}=\A(\omega)\qquad\hbox{a.s in $\omega$.}
\end{equation}
To show that the opposite inclusion holds as well, let  $u$ be the
random function obtained via the Lax--formula \eqref{def u} with
$C_\alpha(\omega)$ in place of $C(\omega)$ and $\overline v$ in
place of $g$. By the minimality property of $r_{\overline v}$,  we
derive that  $r_u(\cdot,\omega)$ is strictly positive on
$C_\alpha(\omega)$ a.s in $\omega$, so combining Corollary
\ref{cor cond C} with Proposition \ref{prop Lax in S} we get that
$u$ is an exact corrector with trace $\overline v(\cdot,\omega)$
on $C_\alpha(\omega)$. Then we invoke Theorem \ref{teo cond C} to
see that $r_{\overline v}(\cdot,\omega)\equiv+ \infty$ on
$C_\alpha(\omega)$ a.s. in $\omega$. This proves that
$C_\alpha(\omega)$ agrees with $\A(\omega)$ a.s. in $\omega$, and
that $\A(\omega)$ is almost surely nonempty. As a consequence we
deduce that $C_\alpha(\omega)$ is almost surely nonempty {\em for
every $\alpha>0$}, see \eqref{eq inclusione}. We can therefore
iterate the above argument to prove that, for every $\alpha>0$,
$C_\alpha(\omega)=\A(\omega)$ a.s. in $\omega$, i.e.
\[
\{ x \in \R^N \,:\, r_{\overline v}(x,\omega) \geq \alpha\}
=
\{ x \in \R^N \,:\, r_{\overline v}(x,\omega)=+\infty\}
\qquad\hbox{a.s. in $\omega$.}
\]
This readily gives
$\R^N\setminus\A(\omega)=\{ x \in \R^N \,:\, r_{\overline v}(x,\omega)=0\}$,
as it was to be shown.

\qed

\bigskip

In the second part of  the section we prove that the random Aubry
set is almost surely foliated by curves defined in $\R$ enjoying
some minimality conditions, where the critical admissible
subsolutions coincide up to an additive constant. When $H$  is
regular enough, such curves turn out to be integral curves of the
Hamiltonian flow. This generalizes properties holding in the
deterministic setting.

\begin{teorema}\label{criti} Assume $\A(\omega) \neq \emptyset$
a.s. in $\omega$. Then there exists a set $\Omega_0$ of
probability 1 such that for any
   $\omega\in\Omega_0$ and any $x\in\A(\omega)$ we can find
a curve $\eta_x:\R\to\A(\omega)$ (depending on $\omega$) with
$\eta_x(0)=x$ satisfying the following properties:
\begin{itemize}
\item[\em (i)] \quad for every $a<b$ in $\R$
\[
S(\eta_x(a),\eta_x(b),\omega)
    =
    \int_a^b (L(\eta_x,\dot\eta_x,\omega) + c)\,\dd s;
\]
\item[\em (ii)] \quad$\displaystyle{\lim_{t \rightarrow \pm \infty} |\eta_x(t)|= + \infty;}$\smallskip
\item[\em (iii)] \quad for every $v\in\S$ the following equality holds a.s. in $\omega$:
\[
 \int_a^b (L(\eta_x,\dot\eta_x,\omega) + c)\,\dd s
    =
    v(\eta_x(b),\omega)-v(\eta_x(a),\omega)\qquad\hbox{for every $a<b$ in $\R$.}\medskip
\]
\end{itemize}
\end{teorema}

We start by some preliminary remarks.

Let $\check H(x,p,\omega) = H(x,-p,\omega)$, and denote by $\check
c$, $\check\S$ and $\check \A(\omega)$ the associated critical
value, the  family of  admissible subsolutions of $\check
H(x,Dv,\omega)=\check c$ and the Aubry set, respectively. It is
easy to see that $\check c=c$. We also have:

\begin{prop}\label{prop A check}
$\check \A(\omega)=\A(\omega)$ a.s. in $\omega$.
\end{prop}

\begin{dimo}
Let $\check S$ the semi--distance associated to $\check H$.
It is easy to check that
\[
\check S(x,y,\omega)=S(y,x,\omega)\qquad\hbox{for every
$x,y\in\R^N$ and $\omega\in\Omega$.}
\]
Let $\overline v$ be a random function of $\S$  weakly strict
outside the Aubry set, see Theorem \ref{teo stretta}. Clearly $-
\overline v \in \check\S$. Let $\Omega_0$ be a set of probability
1 such that for every $\omega\in\Omega_0$ the function $\overline
v(\cdot,\omega)$ is a critical subsolution  and
\[
\R^N\setminus\A(\omega)=\{x\in\R^N\,:\,r_{\overline
v}(x,\omega)=0\,\}.
\]
We claim that the stationary random function $\check r_{-
\overline v}(\cdot,\omega)$, defined through \eqref{def r_v} with
 $-\overline v$ in place of $v$ and $\check S$ in place of $S$,
vanishes in $\R^N \setminus \A(\omega)$ for every
$\omega\in\Omega_0$. This would imply $\check\A(\omega) \subset
\A(\omega)$ a.s. in $\omega$, and arguing analogously  the opposite
inclusion can be obtained as well.

To prove the claim, we argue by contradiction by assuming that
there exist an $\omega\in\Omega_0$ and a point
$x\in\R^N\setminus\A(\omega)$ such that $\check r_{- \overline
v}(x,\omega)>0$. Then there exist an $r>0$ and a point
$y\in\partial B_r(x)$ such that
\begin{equation}\label{eq check}
- \overline v(x,\omega) = - \overline v(y,\omega) + \check
S(y,x,\omega).
\end{equation}
Since $\A(\omega)$ is closed and $-\overline v(\cdot,\omega)$ is a
critical subsolution for $\check H$, we can choose $r>0$ small
enough such that $\partial B_r(x)\subset\R^N\setminus\A(\omega)$.
From \eqref{eq check} we obtain
\[
\overline v(y,\omega)=\overline v(x,\omega)+S(x,y,\omega),
\]
yielding $r_{\overline v}(y,\omega)>0$ with
$y\in\R^\N\setminus\A(\omega)$, a contradiction to the choice
of $\Omega_0$.
\end{dimo}

\bigskip

\ \\
\indent Let
\[ L(x,q,\omega)
    :=
   \max_{p\in\R^N}\left\{\langle p,q\rangle -
   H(x,p,\omega)\right\}. \]
The inequality
   \[
   L(x,q,\omega)\geq
   \max_{H(x,p,\omega)\leq c}\left\{\langle p,q\rangle -
    H(x,p,\omega)\right\}
    =\sigma(x,q,\omega)-c \]
yields
\[
\int_a^b \left(L(\gamma,\dot\gamma,\omega)+c\right)\,\dd t\geq
\int_a^b \sigma(\gamma,\dot\gamma,\omega)\,\dd t
\]
for every curve $\gamma:[a,b]\to\R^N$. We also recall
the relation
\[
S(y,x,\omega)=\inf\left\{\int_0^t
\left(L(\gamma,\dot\gamma,\omega)+c\right)\,\dd
s\,:\,\gamma(0)=x,\,\gamma(t)=y,\,t>0\,\right\}
\]
for every $x,y\in\R^N$ and $\omega\in\Omega$.\smallskip

We approach the proof of Theorem \ref{criti} by proving a weaker
version of it.

\begin{prop}\label{prop prepara}
Assume $\A(\omega) \neq \emptyset$ a.s. in $\omega$ and let $v
\in\S$ be weakly strict outside the Aubry set. Then there exists a
set $\Omega_v$ of probability 1 such that for any
$\omega\in\Omega_v$ and any $x\in\A(\omega)$ we can find a curve
$\eta_x:\R\to\A(\omega)$ (depending on $\omega$) with
$\eta_x(0)=x$ satisfying
\[
  S(\eta_x(a),\eta_x(b),\omega)
    =
    \int_a^b (L(\eta_x,\dot\eta_x,\omega) + c) \,\dd s
    =
v(\eta_x(b),\omega)-v(\eta_x(a),\omega), \] whenever $a<b$ in
$\R$. In addition $\displaystyle{\lim_{t \rightarrow \pm \infty}
|\eta_x(t)|= + \infty}$.

\end{prop}

\begin{dimo}
We take $\Omega_v$ such that for every $\omega\in\Omega_v$ the
function $v(\cdot,\omega)$ is a critical subsolution,
$\A(\omega)\not =\emptyset$ and
\[
\A(\omega)= \{x\in\R^N\,:\,r_v(x, \omega) = + \infty\},\quad\R^N\setminus
\A(\omega) = \{x\in\R^N\,:\,r_v(x, \omega) = 0\}.
\]
Fix $\omega\in\Omega_v$. The function
\[
u(x):=\inf\{v(x,\omega)+S(y,x,\omega)\,:\,y\in \A(\omega)\,\},
\qquad\hbox{$x\in\R^N$,}
\] is a
viscosity solution of
\[
H(x,Du,\omega)=c\quad\hbox{in $\R^N$,}
\]
and consequently $u(x) - c \,t$ is  a solution of the
time--dependent Hamilton--Jacobi Cauchy problem
\begin{eqnarray*}
\begin{cases}
\partial_t w + H(x,D w,\omega)=0 & \hbox{in
$(0,+\infty)\times\R^N$}\\
w(0,x,\omega)=u(x) & \hbox{in $\R^N$}.
\end{cases}
\end{eqnarray*}
Hence the following Lax--Oleinik representation formula holds for
every $x\in\R^N$ and $t>0$:
\begin{equation}\label{lax-oleinik formula}
u(x)=\inf\left\{ u(\gamma(-t))+\int_{-t}^0
(L(\gamma(s),\dot\gamma(s),\omega) + c)\,\dd
s\,:\,\gamma(0)=x\,\right\},
\end{equation}
where $\gamma$ varies in the family of absolutely continuous
curves from $[-t,0]$ to $\R^N$, see \cite{DavSic05}. By standard
arguments of the Calculus of Variations \cite{BuGiHi98}, a
minimizing absolutely continuous curve does exist for any fixed $t
> 0$ thanks to the coercivity and lower semicontinuity  properties
of $L$. Moreover such curves turn out to be equi--Lipschitz
continuous, see \cite{D1-06}. Given an increasing sequence  $t_n$
with $\lim_n t_n=+\infty$, we denote by $\gamma_n$ the
corresponding minimizers and extend them on the whole interval
$(-\infty,0]$ by setting  $\gamma_n(t)= \gamma_n(-t_n)$ in $(-
\infty, -t_n)$, for any $n$.  Thanks to Ascoli Theorem, the
sequence $\gamma_n$ so defined  has a local uniform limit, denoted
by $\gamma_x$, in $( - \infty,0]$, up to subsequences.  Taking
into account the optimality of the $\gamma_n$ and the fact that
$u$ is a critical (sub)solution, we get for any $t>0$
\begin{equation}\label{curva calibrata}
u(\gamma_x(0))-u(\gamma_x(-t))=\int^0_{-t}
(L(\gamma_x,\dot\gamma_x,\omega) + c)\,\dd s
=S(\gamma_x(-t),\gamma_x(0),\omega),
\end{equation}
and for any $a<b\leq 0$,
\[
u(\gamma_x(b))-u(\gamma_x(a))=\int_a^{b}
(L(\gamma_x,\dot\gamma_x,\omega) + c)\,\dd s
=S(\gamma_x(a),\gamma_x(b),\omega).
\]
If, in particular, $x = \gamma_x(0)\in\A(\omega)$,  we have
\[
u(\gamma_x(0))=v(\gamma_x(0),\omega),\qquad u(\gamma_x(-t))\geq
v(\gamma_x(-t),\omega)\quad\hbox{for every $t>0$.}
\]
From \eqref{curva calibrata} we then derive
\[
S(\gamma_x(-t),\gamma_x(0),\omega)
    \leq
    v(\gamma_x(0),\omega)-v(\gamma_x(-t),\omega),
\]
which in turn implies that $v(\cdot,\omega)$ and $u(\cdot)$
coincide on $\gamma_x$. Since $r_v(\cdot,\omega)$ vanishes outside
$\A(\omega)$, we conclude  the support of $\gamma_x$ is contained
in $\A(\omega)$, as claimed.\par

The same argument  can be applied to the function
$-v(\cdot,\omega)$ and to the Hamiltonian $\check
H(x,p,\omega):=H(x,-p,\omega)$. In view of Proposition \ref{prop A
check}, we can assume, without any loss of generality,  that $\check
\A(\omega)=\A(\omega)$.  Taking into account the relations
\[
\check L(x,q,\omega)=L(x,-q,\omega),\quad
\check\sigma(x,q,\omega)=\sigma(x,-q,\omega),\quad\check
S(x,y,\omega)=S(y,x,\omega)
\]
for every $x,\,y,\,q\in\R^N$,  we deduce as above that for  every
$x\in\A(\omega)$ there exists a curve
$\xi_x:(-\infty,0]\to\A(\omega)$ with $\xi_x(0)=x$ satisfying
\[
-v(\xi_x(b),\omega)+v(\xi_x(a),\omega)=\int^b_{a} \check
L(\xi_x,\dot\xi_x,\omega)\,\dd s=\check
S(\xi_x(a),\xi_x(b),\omega)
\]
for every $a<b\leq 0$. The curve $\eta_x$ with the claimed
properties  is obtained by setting
\[
\eta_x(t):=
\begin{cases}
\xi_x(-t) & \hbox{if $t\geq 0$}\\
\gamma_x(t) & \hbox{if $t\leq 0$.}
\end{cases}
\]
Finally, if there is a limit point $\widetilde x$ of
$\eta_x$ for $t \rightarrow \pm \infty$,  we can find a sequence
of compact intervals $[a_n,b_n]$ with $b_n-a_n\geq n$ such that $\eta_x(a_n)$ and $
\eta_x(b_n)$ both converge to $\widetilde x$ and
\[
\int_{a_n}^{b_n} (L(\eta_x,\dot\eta_x,\omega) + c)\,\dd s
=S(\eta_x(a_n),\eta_x(b_n),\omega) \rightarrow 0.
\]
By joining $\widetilde x$ to $\eta_x(a_n)$ and to $\eta_x(b_n)$ with two segments,
we can define a sequence of loops $\xi_n:[0,t_n]\to\R^N$ with $\widetilde x$ as base
point such that $t_n\to +\infty$ and
\[
\int_{0}^{t_n} (L(\xi_n,\dot\xi_n,\omega) + c)\,\dd s\to 0.
\]
This would imply that $\widetilde x \in \A_f(\omega)$ in view of
\eqref{aub}. Since $\A_f(\omega)$ is almost surely empty by
hypothesis, we see that no such points can exist and the limit
relation at infinity asserted in the statement follows.
\end{dimo}

\ \\
\indent We proceed to show that the minimal curves  $\eta_x$
 can be chosen independently of $v\in\S$.

\smallskip

{\bf\noindent Proof of Theorem \ref{criti}.}
 The critical subsolution $\overline v$ appearing
 in the statement of Proposition \ref{stramin} is weakly strict outside the Aubry set and  has the form
 \[
\overline v(x,\omega)=\sum_n \lambda_n
v_n(x,\omega)\qquad\hbox{for every
$(x,\omega)\in\R^N\times\Omega$,}
\]
where the $(\lambda_n)_n$ are positive constants  satisfying
$\sum_n \lambda_n=1$, and $(v_n)_n$ is a sequence dense in
$\widehat\S$ with respect to the the almost sure convergence in
$\D C(\R^N)$. By
Proposition \ref{prop prepara}  there exists a set $\Omega_0$ of
probability 1 such that, for every $\omega\in\Omega_0$ and every
$x\in\A(\omega)$, we can find a curve $\eta_x:\R\to\R$ satisfying
$\eta_x(0)=x$, $\lim_{t \rightarrow \pm \infty} |\eta_x(t)|= +
\infty$ and
\begin{equation}\label{cond critica}
   \overline  v(\eta_x(b),\omega)- \overline v(\eta_x(a),\omega)
    =
    \int_a^b (L(\eta_x,\dot\eta_x,\omega) + c)\,\dd s
    =
    S(\eta_x(a),\eta_x(b),\omega)
\end{equation}
whenever $a<b$ in $\R$. Since $v_n\in \S$, up to removing from $\Omega_0$ a set
of probability 0, we can furthermore assume that, for any $\omega\in\Omega_0$, each
function $v_n(\cdot,\omega)$ is a subsolution of \eqref{eq critica}. This readily
implies that, for any such $\omega$, equality \eqref{cond critica} holds with
$v_n$ in place of $\overline v$, for every $n\in\N$.

To prove {\em (iii)}, fix $v\in\S$ and set $\hat v=v-v(0,\omega)$.
Clearly, it suffices to show the assertion for $\hat v$. Since
$\hat v\in\widehat \S$, there exists a subsequence $(v_{n_k})_k$
and a set $\widehat\Omega\subseteq\Omega_0$ of probability 1 such
that $v_{n_k}(\cdot,\omega)\ucv v(\cdot,\omega)$ for any
$\omega\in\widehat\Omega$. By passing to the limit, we derive that
equality \eqref{cond critica} holds with $\hat v$ in place of
$\overline v$ for any such $\omega\in\widehat\Omega$, as it was to
be shown.
 \qed
\medskip
\end{section}

\begin{section}{Open questions}\label{sez open}

This is the third of a series of papers we have devoted to the
analysis of critical equations for stationary ergodic
Hamiltonians, see  \cite{DavSic08, DS1-08}, by using the metric
approach combined with some tools from Random Set Theory. This
method has allowed to get a complete picture of the setup when the
state variable space is $1$--dimensional, as specified in the
introduction,  and, we think, has revealed to be effective also in
the multidimensional setting, highlighting some interesting
analogies with the compact case. However many crucial problems are
still to be clarified. The more striking is:

\medskip

(1) {\em In case of existence of an exact corrector, is the random
Aubry set almost surely nonempty ?}

\medskip

In view of Theorem \ref{teo stretta}, we can put it more
dramatically:

\medskip

(1$'$) { \em Is it impossible the simultaneous existence of an
exact corrector and a global weakly strict admissible critical
subsolution ?}

\medskip

In this respect, it should be helpful to strengthen Theorem
\ref{teo stretta}, as in periodic case. So we would  also like to
know:

\medskip

(2) {\em If the Aubry set is a.s. empty, there exist {\em strict}
global critical subsolutions?  Can we find one of such subsolution
which is, in addition, smooth?}

\medskip

If the answer to (1), (1$'$) is positive, another question urges
itself upon us:

\medskip

(3) {\em Is any exact corrector the Lax extension from the Aubry
set of an admissible trace ?}

\medskip
Or, in other terms, is the Aubry set an uniqueness set for the
critical equation, as in the deterministic compact case? Notice
that both questions (1) and (3) have positive answer when $N=1$,
see \cite{DavSic08}, and in any space dimension when $c=c_f=
\sup_x\,\min_p H(x,p,\omega)$
    and the critical stable norm in nondegenerate, see \cite{DS1-08}.\par

On the contrary, if $c=c_f$, $\A_f(\omega)$ is a.s. empty and the
critical stable norm is nondegenerate, then no exact solutions can
exist. We stress that, as far as we know, all counterexamples
published in the literature to the existence of exact correctors
are in this frame. It should be interesting to find, if possible,
counterexamples in cases where the previous conditions are not
satisfied.\par

 The above  nonexistence result morally says that
we can hope to have exact correctors only  if some metric
degeneracy of $S_c$ takes place either at finite points (i.e. when
$\A_f(\omega) \neq \emptyset$ a.s. in $\omega$) or at infinity
(i.e. when the stable norm vanishes in some directions).\par

 The converse is partially true, in the sense
that when $c=c_f$ and the classical Aubry set $\A_f(\omega)$ is
almost surely nonempty, we know that exact correctors do exist.
One should wonder if a corrector does exist in case of sole metric
degeneracy at infinity. Indeed, this is unclear even if $c=c_f$
and, evidently, $\A_f(\omega)$ is a.s. empty. We have exhibited an
example in \cite{DS1-08}, see Example 6.8, of Eikonal equation of
the type
\begin{equation}\label{eq Eikonal}
|Du(x,\omega)|^2=V(x,\omega) \qquad\hbox{in $\R^N$,}
\end{equation}
where the potential $V$ is a random continuous stationary bounded
positive function with  infimum a.s. equal to $0$, with the
peculiarity that the corresponding critical stable norm is equal
to 0, i.e. vanishes in any direction. Note that here $c=c_f=0$
since the null function is a strict admissible subsolution, and no
subsolutions of \eqref{eq HJa} exist for $a<0$. Here we face a
dilemma: either an exact corrector does exist, and then the
question (1), (1$'$) has a negative answer, since the Aubry set
$\A(\omega)$ must be a.s. empty (for the null function is a strict
admissible subsolution); or we have to recognize that metric
degeneracy at infinity is not sufficient for yielding critical
solutions.\par

We remark that a negative answer to questions (1), (1$'$) would
come from the finiteness of
\[
\liminf_{|y|\to +\infty}S(y,0,\omega) \quad\text{ a.s. in
$\omega$},
\]
where $S$, as usual, is the critical distance associate to
\eqref{eq Eikonal}. In fact, if such limit is less than $+
\infty$, then, by the triangle inequality and  other properties
enjoyed by $S$, it is easily seen that
\[
\Omega_0:=\{\omega\,:\,\liminf_{|y|\to
+\infty}S(y,x,\omega)<+\infty\quad\hbox{for some $x\in\R^N$}\},
\]
has probability $1$, and so a finite--valued random function $u$
can be defined by setting
\[
u(x,\omega)=\liminf_{|y|\to +\infty}S(y,x,\omega)\qquad\hbox{for
$\omega\in\Omega_0$}
\]
and $u(x,\omega)=0$ otherwise, for every $x\in\R^N$.  Via standard
arguments, it can be then proved that $u(\cdot,\omega)$ is a
solution of \eqref{eq Eikonal}.\par

Another subject of interest is about approximate correctors. So
far we don't have any counterexamples to their existence when
exact correctors do not exist. Hence the main question is:
\medskip

(4) {\em Do approximate correctors always exist?}

\medskip

This issue  is also strongly related to homogenization problems
and a positive answer would be an important step towards
generalizations of the results proved in \cite{ReTa00, Souga99} to
more general Hamiltonians.\par

As usual, the answer is positive if $N =1$, or in any space
dimension if  $c=c_f= \sup_x\,\min_p H(x,p,\omega)$. In this
setting, we have in addition proved that approximate correctors
can be represented by Lax formula \eqref{def u}, taking as random
source set the $\delta$--maximizers over $\R^N$ of the function $x
\mapsto \min_p H(x,p,\omega)$. This result essentially exploits
the assumption that $c_f$ is the supremum of such function, which
is always the case in the $1$--dimensional setting.\par

To extend it in more general setups, at least when $c=c_f$, the
idea could be to replace the $\delta$--minimizers by some sort of
{\em approximate} Aubry set. But such a set seems not easy to
define even assuming the existence of a smooth strict critical
subsolution,  and  so this attempt has not given, till now, any
output.
\par

Note that the existence results of \cite{Is}  for approximate
correctors in the almost--periodic case are based on an ergodic
approximation of the Hamilton--Jacobi equation, and so are not
constructive. A  final question, which stems from the previous
discussion,  then is

\medskip

(5) {\em  At least in the almost--periodic case, are the
approximate correctors representable through Lax formulae?}

\medskip

\end{section}

\begin{appendix}
\begin{section}{}

We begin recalling that a function $f$ defined on $\R^N$ is said
to be {\em almost--periodic} if it is bounded, continuous and if
it can be approximated, uniformly on $\R^N$, by finite linear
combinations of functions in the set $\{\e^{2\pi i
\langle\lambda,\,x\rangle}\,:\,\lambda\in\R^N\,\}$, see
\cite{AmFrid,Ru90} for instance.

\smallskip

This appendix is devoted to show that any  almost--periodic
Hamiltonian is a specific realization of a stationary ergodic
Hamiltonian, with underlying probability space $\Omega$  {\em
separable}  in a measure theoretic   sense. Generalizing the
construction of the quasi--periodic case, we more precisely  prove
that $\Omega$ can be taken as the infinite--dimensional  torus
with $\R^N$ appropriately acting on it. Therefore  $\Omega$ is in
addition   a {\em compact metric space} with the product topology,
and is as well separable
 from a  topological viewpoint.

 \smallskip

The  statement of the main result  is the following:

\begin{teorema}\label{teo AP realization}
Let $H_0:\R^N\times\R^N\to\R$ be a continuous Hamiltonian
satisfying the following assumptions:
\begin{itemize}
   \item[(B1)] \ $H_0(\cdot,p)$ is almost--periodic in $\R^N$ for every fixed
     $p\in\R^N$;\smallskip
      \item[(B2)]
    $\ H_0(x,\cdot)\ \hbox{is  convex on $\R^N$ for every $x\in\R^N$;}$ \smallskip
     \item[(B3)] \ there
     exist two superlinear continuous functions $\alpha,\beta:\R_+\to\R$ such
     that
     \[
     \alpha\left(|p|\right)\leq H_0(x,p)\leq \beta\left(|p|\right)\qquad\hbox{for all
     $(x,p)\in\R^N\times\R^N$;}
     \]
      \item[(B4)] the set of minimizers of $H_0(x,\cdot)$
    has empty interior for every $x\in\R^N$.\smallskip
\end{itemize}
Then there exist a {\em separable} probability space
$(\Omega,\F,\PP)$, an ergodic group of translations
$(\tau_x)_{x\in\R^N}$ and an Hamiltonian
$H:\R^N\times\R^N\times\Omega\to\R$ satisfying assumptions
(H1)--(H5) of Section \ref{sez HJ} such that
\begin{equation*}
H(x,p,\omega_0)=H_0(x,p)\qquad\hbox{for every
$(x,p)\in\R^N\times\R^N$,}
\end{equation*}
for some $\omega_0\in\Omega$.
\end{teorema}

\indent The proof of Theorem \ref{teo AP realization} will require
some  preliminary work. We start by  some classical definitions
and results from the theory of dynamical systems, see \cite{KaHa}.
\smallskip\\
\indent A continuous map $\tau:\Omega\to\Omega$ defined on a {\em
Hausdorff topological space} $\Omega$ will be said to be {\em
minimal} if the orbit
$$\D{orb}(\omega):=\{\tau^n(\omega)\,:\,n\in\Z\}$$
of every point $\omega\in\Omega$ is dense in $\Omega$.
%
A Borel probability measure $\mu$ on $\Omega$ is called {\em
$\tau$--invariant} if $\mu(\tau^{-1}(E))=\mu(E)$ for every
$\mu$--measurable set $E$. A measurable subset $E$ of $\Omega$ is
called $\tau$--{\em invariant} if $\mu(\tau^{-1}(E)\,\Delta\,
E)=0$, where $\Delta$ stands for the symmetric difference. A
$\tau$--invariant measure $\mu$ is called {\em ergodic} (with
respect to $\tau$) if for any $\tau$--invariant measurable set
$E\subset \Omega$ either $\mu(E)=0$ or $\mu(E)=1$. When $\Omega$
is a metrizable compact space, $\tau$ will be called {\em uniquely
ergodic} if it has only one invariant Borel probability measure
$\mu$. In this instance $\mu$ is necessarily ergodic with respect
to $\tau$, see \cite[Proposition 4.1.8]{KaHa}.

We will use the following result from \cite[Proposition
4.1.15]{KaHa}.

\begin{prop}\label{a prop uniquely ergodic}
Let $\Omega$ be a metrizable compact space and $\tau : \Omega
\rightarrow \Omega$ a continuous map. If for every continuous
function $\varphi$ belonging to  a dense set  in the space
$\D{C}(\Omega)$ the time averages
$(1/n)\sum_{k=0}^{n-1}\varphi(\tau^k(\omega))$ converge uniformly
to a constant, then $\tau$ is uniquely ergodic.
\end{prop}

By applying  Proposition \ref{a prop uniquely ergodic}, we show

\begin{prop}\label{a prop unique erg}
Let $\Omega$ and $\tau$  be a compact metric space and an isometry
on it, respectively. If $\tau$ has a dense orbit, then it is
uniquely ergodic.
\end{prop}

\begin{dimo}
Let $\varphi\in \D{C}(\Omega)$. In view of Proposition \ref{a prop
uniquely ergodic} it suffices to show that the functions
\[
\varphi_n(\omega)=\frac{1}{n}\sum_{k=0}^{n-1}\varphi(\tau^k(\omega))
\]
uniformly converge to a constant. It is easy to see that the
functions $\varphi_n$ are equi--bounded by $\|\varphi\|_\infty$,
which is finite since $\Omega$ is compact and $\varphi$ is
continuous. Moreover, they are equi--continuous, because a
continuity modulus for $\varphi$  plays the same role for each of
the $\varphi_n$, since $\tau$ preserves the distance. By
Ascoli--Arzel\`a Theorem we infer that $\varphi_n$ uniformly
converge to a function $\psi$ which is $\tau$--invariant, i.e.
constant on the orbits of $f$. Since there is a dense orbit by
hypothesis and $\psi$ is continuous, we conclude that $\psi$ is
constant, as it was to be proved.
\end{dimo}

\medskip

Note that the previous result applies, in particular, to minimal
maps, for which  all the orbits are dense.

\bigskip

Let $\T^1$ be the one--dimensional flat torus  endowed with the
flat Riemannian metric, still denoted by $|\cdot|$,  induced by
the Euclidean metric on $\R$. We define a distance $d$ on
$\T^\infty:=\Pi_{j=1}^{+\infty}\T^1$  via
\begin{equation}\label{a def dist}
d(\omega,\omega')=\sum_{n=1}^{+\infty}\frac{1}{2^n}|\omega_n-\omega'_n|
\qquad\hbox{$\omega=(\omega_n)_n$,\, $\omega'=(\omega'_n)_n$ in
$\T^\infty$.}
\end{equation}
 By Tychonoff Theorem, $\T^\infty$ is a compact metric space with
respect to $d$.  We consider  on $\T^\infty$ the product
probability measure $\mu:=\Pi_{j=1}^{+\infty}\leb^1\restr\T^1$.
For every $m\in\N$ we denote by $\pi_m:\T^\infty\to\T^m$ the
projection on the first $m$--components, and   by
$\mu_m:={\pi_m}_\natural\,\mu$ the push--forward on $\T^m$ of the
measure $\mu$, i.e. the probability measure given by
\[
\mu_m(E)=\mu\left(\pi_m^{-1}(E)\right)\qquad \hbox{for every Borel
set $E\subseteq\T^m$.}
\]
We endow $\T^m$ with the distance $d_m$ defined as
\[
d_m(\omega,\omega')=\sum_{j=1}^{m}\frac{1}{2^j}|\omega_j-\omega'_j|,\qquad\hbox{$\omega$,\,
$\omega'\in\T^m$.}
\]

Given a sequence $(\lambda_n)_n$ of vectors in $\R^N$, we consider
the group  of translation $(\tau_x)_{x\in\R^N}$ defined as
\begin{equation}\label{a eq translations}
(\tau_x\omega)_j\equiv
\omega_j+\langle\lambda_j,x\rangle\quad\hbox{(mod
1)}\qquad\hbox{for every $j\in\N$,}
\end{equation}
 Note
that $\mu$ is invariant with respect to $(\tau_x)_{x\in\R^N}$. We
denote, for any  $x\in\R^N$, by ${\tau_x}_{|\T^m}:\T^m\to\T^m$ the
translation $\tau_x$ restricted to the first $m$ components, i.e.
\[
({\tau_x}_{|\T^m}\,\omega)_j\equiv
\omega_j+\langle\lambda_j,x\rangle\quad\hbox{(mod
1)}\qquad\hbox{for every $j\in\{1,\dots,m\}$,}
\]
for each $\omega=(\omega_1,\dots,\omega_m)\in\T^m$. Clearly
$\mu_m$ is invariant with  respect to
$({\tau_x}_{|\T^m})_{x\in\R^N}$.

Motivated by the next result, we are specially interested to the
case where the sequence  $(\lambda_n)_n$ in $\R^N$ is {\em
rationally independent}, i.e. when every finite combination of
elements of the sequence with rational coefficients is zero if and
only all the coefficients vanish.
\smallskip\\
\indent The following holds

\begin{prop}\label{a prop translation+}
Let $(\lambda_n)_n$ be a countable family of rationally
independent vectors in $\R^N$. Then there exists $\hat x\in\R^N$
such that the translations ${\tau_{\hat{x}}}$ and
${\tau_{\hat{x}}}_{|\T^m}$ are minimal on $\T^\infty$ and on
$\T^m$ for every $m\in\N$, respectively. In particular, $\mu$ and
$\mu_m$ are uniquely ergodic with respect to ${\tau_{\hat x}}$ and
${\tau_{\hat x}}_{|\T^m}$, respectively.
\end{prop}

We will exploit in  the proof the following known fact, see
\cite[Proposition 1.4.1]{KaHa}.

\begin{prop}\label{a prop translation}
Let $\gamma=(\gamma_1,\dots,\gamma_m)$ be a vector of $\R^m$ and
let $T_\gamma$ be the translation on the torus $\T^m$ defined as
\[
T_\gamma(\omega_1,\dots,\omega_m)\equiv(\omega_1+\gamma_1,\dots,\omega_m+\gamma_m)\quad\hbox{\rm
(mod 1).}
\]
Then $\T_\gamma$ is minimal if and only if $\sum_{j=1}^m k_j
\gamma_j\not\in\Z$ for any choice of $(k_1,\dots,k_m)$ in
$\Z^m\setminus\{ 0\}$.
\end{prop}

\noindent{\bf Proof of Proposition \ref{a prop translation+}.} Let
us consider the countable set
\[
\mathcal I:=\{\,\Bbbk=(k_n)_n\in\Z^\N\,:\, k_j\not=0\ \hbox{for a
finite and positive number of indices $j$\,}\}.
\]
For every $\Bbbk\in\mathcal I$,  we define
\[
V_\Bbbk:=\{x\in\R^N\,:\,\sum_i k_i \,\langle
\lambda_i,x\rangle\not\in\Z\,\}
\]
since $\sum_i k_i\lambda_i\not=0$, this set is open and dense in
$\R^N$. Baire's Theorem then implies that
$V:=\cap_{\Bbbk\in\mathcal I}\, V_\Bbbk$ is dense, in particular
is non void. Pick $\hat x\in V$. The minimality of ${\tau_{\hat
x}}_{|\T^m}$ in $\T^m$ for every $m\in\N$ follows from Proposition
\ref{a prop translation}.

Let us show that $\tau_{\hat x}$ is minimal in $\T^\infty$, i.e.
\[
\D{orb}(\omega)\cap B_r(\omega')\not=\emptyset \quad\text{for any
$\omega$, $\omega'$ in $\T^\infty$, any $r >0$.}
\]
 Let $m\in\N$ be large enough to have
$\sum_{j=m+1}^\infty 1/2^j<r/2$. Since ${\tau_{\hat x}}_{|\T^m}$
is minimal on $\T^m$, there exists an integer $k\in\Z$ such that
\[
\sum_{j=1}^{m}\frac{1}{2^j}\,\|(\tau_{\hat
x}^k(\omega))_j-\omega'_j\|<r/2.
\]
Hence
\[
d(\tau_{\hat x}^k(\omega),\omega')
    =
\sum_{j=1}^{\infty}\frac{1}{2^j}\,\|(\tau_{\hat
x}^k(\omega))_j-\omega'_j\|
    \leq
\sum_{j=1}^{m}\frac{1}{2^j}\,\|(\tau_{\hat
x}^k(\omega))_j-\omega'_j\|
    +
\sum_{j=m+1}^{\infty}\frac{1}{2^j}
    <
    r.
\]
The remainder of the assertion is a straightforward consequence of
Proposition \ref{a prop unique erg}. \qed

\medskip
We summarize what we have proved so far in the next statement.

\begin{teorema}
Let $(\lambda_n)_n$ be a countable family of rationally
independent vectors in $\R^N$, $d$ the distance on $\T^\infty$
defined via \eqref{a def dist}, $(\tau_x)_{x\in\R^N}$ the group of
translations on $\T^\infty$ defined according to \eqref{a eq
translations},  and $\mu$ the product probability measure defined
as $\mu:=\Pi_{j=1}^{+\infty}\leb^1\restr\T^1$. Then
$(\T^\infty,d)$ is a compact metric space, in particular
separable, and $(\tau_x)_{x\in\R^N}$ is ergodic with respect to
$\mu$.
\end{teorema}

We proceed to show that given  an almost--periodic function $f$ on
$\R^N$, a  sequence of rationally  independent vectors
$(\lambda_n)_n$  can be chosen in such a way that  $f$ is a
specific realization of a random variable  on $\T^\infty$ with
respect to the group of translations $(\tau_x)_{x\in\R^N}$ defined
via \eqref{a eq translations}; in addition such random variable
can be taken continuous. In the sequel, we will denote by $\mathbf
0$ the element of $\T^\infty$ all of whose components are equal to
0.

\begin{prop}\label{a prop approximation}
Let $f$ be an almost periodic function in $\R^N$. There exist a
sequence of rationally independent vectors $(\lambda_n)_n$ in
$\R^N$,   inducing a dynamical system $(\tau_x)_{x\in\R^N}$ on
$\T^\infty$ via \eqref{a eq translations}, and  a continuous
function $\underline f:\T^\infty\to\R$ such that $f(x)=\underline
f(\tau_x\mathbf 0)$.
\end{prop}

\begin{dimo}
In what follows, we will use some known facts about
almost--periodic functions, see \cite{Ru90}. For every
$\lambda\in\R^N$, let us set
\[
a_\lambda:=\lim_{R\to +\infty}\int_{B_R}\!\tagli f(x)\,\e^{-2\pi
i\langle\lambda,x\rangle}\,\dd x
\]
and $\Lambda:=\{\,\lambda\in\R^N\,:\,a_\lambda\not=0\,\}$. Since
$f$ is almost periodic,  the set $\Lambda$ is countable,  and  we
will write $\Lambda=(\widetilde\lambda_n)_n$. For every $n\in\N$,
we define
\[
f_n(x)=\sum_{k=1}^n a_{\widetilde\lambda_k} \e^{2\pi
i\langle\widetilde\lambda_k,x\rangle}
\]
It is well  known that $f_n$ converge uniformly to $f$ in $\R^N$.
We now want to write $f$ as limit of a totally convergent series.
To this purpose, we choose an increasing sequence of integers
$(\widetilde k_n)_n$ in such a way that $\|f_{\widetilde
k_n}-f\|_{L^\infty(\R^N)}\leq 1/2^{n+2}$, and we set
$g_1(\cdot)=f_{\widetilde k_1}(\cdot)$ and, for $n\geq 2$,
\[
g_n(x)
    :=
    f_{\widetilde k_n}(x)-f_{\widetilde k_{n-1}}(x)
    =
\sum_{j=\widetilde k_{n-1}+1}^{\widetilde k_n}
a_{\widetilde\lambda_j} \e^{2\pi
i\langle\widetilde\lambda_j,x\rangle}\qquad x\in\R^N.
\]
Clearly $f(x)=\sum_{n=1}^\infty g_n(x)$. Furthermore,
$\|g_n\|_{L^\infty(\R^N)}\leq C/2^n$ for every $n\in\N$, where $C$
is a constant greater than $1+2\|f_{\widetilde
k_1}\|_{L^\infty(\R^N)}$.  From $(\widetilde\lambda_n)_n$ we
extract a sequence $(\lambda_n)_n$ of vectors rationally
independent in such a way that each $\widetilde\lambda_n$ is a
rational linear combination of $\lambda_1,\dots,\lambda_n$. By
expressing each $\widetilde\lambda_j$ in $g_n$ in terms of its
rational finite linear combination via elements of
$(\lambda_n)_n$, we derive that
\[
g_n(x)=G_n\left(\langle\lambda_1,x\rangle,\dots,\langle\lambda_{k_n},x\rangle\right),\quad
x\in\R^N,
\]
where $(k_n)_n$ is a non decreasing sequence of indexes with
$k_n\leq \widetilde k_n$, and $G_n(\omega_1,\dots,\omega_{k_n})$
is a continuous function from $\T^{k_n}$ to $\C$. For every $n$,
we define a continuous function $\underline g_n$ on $\T^\infty$ by
setting
\[
\underline
g_n(\omega)=G_n\comp\pi_{k_n}(\omega),\qquad\hbox{$\omega\in\T^\infty$.}
\]
Let $(\tau_x)_{x\in\R^N}$ be the group of translations on
$\T^\infty$ associated with the vectors $(\lambda_n)_n$ via
\eqref{a eq translations}. Note that $\underline g_n(\tau_x\mathbf
0)=g_n(x)$ for every $x\in\R^N$. Since $\{\,\tau_x(\mathbf
0)\,:\,x\in\R^N\,\}$ is dense in $\T^\infty$ by Proposition \ref{a
prop translation+} and $\underline g_n$ is continuous on
$\T^\infty$, we derive that $\|\underline
g_n\|_{L^\infty(\T^\infty)}\leq C/2^n$. This yields that the
series
\[
\sum_{n=1}^{+\infty} \underline
g_n(\omega),\qquad\omega\in\T^\infty
\]
uniformly converges to a continuous function $\underline
f:\T^\infty\to\C$, in particular
\[
    \underline f(\tau_x\mathbf 0)
    =
    \sum_{n=1}^{+\infty} \underline g_n(\tau_x\mathbf 0)
    =
    \sum_{n=1}^\infty g_n(x)
    =
    f(x)\qquad \hbox{for every $x\in\R^N$.}
\]
The fact that $\underline f(\T^\infty)\subset\R$ finally follows
by noticing that the continuous function $\underline f$ takes real
values on $\{\,\tau_x(\mathbf 0)\,:\,x\in\R^N\,\}$, which is dense
in $\T^\infty$.
\end{dimo}
\ \\

\indent The last step consists in extending the previous result to
functions that additionally depend on $p$.

\begin{prop}\label{a prop preliminary}
Let $H_0:\R^N\times\R^N\to\R$ be a continuous function satisfying
the following assumptions:
\begin{itemize}
    \item[(A1)] for every $p\in\R^N$ the function $H_0(\cdot,p)$ is
almost--periodic in $\R^N$;\smallskip
    \item[(A2)] for every $R>0$ there exists a modulus $\eta_R$
    such that
    \[
    \quad |H_0(x,p)-H_0(x,q)|\leq \eta_R(|p-q|)\qquad\hbox{for every $x\in\R^N$ and $p,q\in B_R$.}
    \]
\end{itemize}
Then there exists a continuous $\underline
H:\T^\infty\times\R^N\to\R$ such that
\begin{equation*}
\underline H(\tau_x\mathbf 0,p)=H_0(x,p)\qquad\hbox{for every
$(x,p)\in\R^N\times\R^N$,}
\end{equation*}
where $(\tau_x)_{x\in\R^N}$ is the group of translations on
$\T^\infty$ defined according to \eqref{a eq translations} for a
suitable chosen sequence of rationally independent vectors
$(\lambda_n)_n$ in $\R^N$.
\end{prop}

\begin{dimo}
For every $\lambda$ and $p$ in $\R^N$ let us set
\[
a_\lambda(p):=\lim_{R\to +\infty}\int_{B_R}\!\tagli
\,H_0(x,p)\,\e^{-2\pi i\langle\lambda,x\rangle}\,\dd x.
\]
The fact that $a_\lambda$ is continuous on $\R^N$ for every fixed
$\lambda$,   follows from the estimate
\[
|a_\lambda(p)-a_\lambda(q)|\leq\eta_R(|p-q|),
\]
which holds for every  $R>0$, $p,q\in B_R$. Let $(p_k)_k$ be a
dense sequence in $\R^N$ and set
\[
\Lambda:=\bigcup_{k\in\N}\{\,\lambda\in\R^N\,:\,a_\lambda(p_k)\not=0\,\}.
\]
The almost--periodic character of $H_0(\cdot,p)$ implies that
$\Lambda$ is countable, so we will write
$\Lambda=(\widetilde\lambda_n)_n$. From the continuity of
$a_\lambda$ we deduce
\[
a_\lambda(\cdot)\equiv 0\qquad\hbox{for every
$\lambda\not\in\Lambda$.}
\]
From $(\widetilde\lambda_n)_n$ we extract a sequence
$(\lambda_n)_n$ of rationally independent  vectors  in such a way
that each $\widetilde\lambda_n$ is a rational linear combination
of $\lambda_1,\dots,\lambda_n$. Let $(\tau_x)_{x\in\R^N}$ be the
group of translations on $\T^\infty$ associated to the vectors
$(\lambda_n)_n$ via \eqref{a eq translations}. In view of
Proposition \ref{a prop approximation}, for every $p\in\R^N$ there
exists a continuous function $\underline
H(\cdot,p):\T^\infty\to\R$ such that
\[
\underline H(\tau_x\mathbf 0,p)=H_0(x,p)\qquad\hbox{for every
$x\in\R^N$.}
\]
From this we get that, for every $\omega\in\{\,\tau_x(\mathbf
0)\,:\,x\in\R^N\,\}$,
\[
|\underline H(\omega,p)-\underline H(\omega,q)|\leq
\eta_R(|p-q|)\qquad\hbox{for every $p,q\in B_R$ and $R>0$.}
\]
Since $\{\,\tau_x(\mathbf 0)\,:\,x\in\R^N\,\}$ is dense in
$\T^\infty$ and $\underline H(\cdot, p)$ is continuous on
$\T^\infty$ for every fixed $p$, we derive that the above
inequality holds for every $\omega\in\T^\infty$. Hence $\underline
H$ is jointly continuous in $(\omega,p)$ and the proof is
complete.
\end{dimo}

\medskip
We are now in position to prove Therem \ref{teo AP realization}.\medskip\\
\noindent{\bf Proof of Theorem \ref{teo AP realization}.} We
recall (see for instance \cite{R70}) that a convex function
$\psi:\R^N\to\R$ is locally Lipschitz, and its Lipschitz constant
in $B_R$ can be controlled with the supremum of $|\psi|$ on
$B_{R+2}$, for every $R>0$. In particular the Hamiltonian $H$
satisfies assumption (A2) in Proposition \ref{a prop preliminary}
with $\eta_R(h):=L_R\,h$, where
\[
L_R:=\sup\{\,|H_0(x,p)|\,:\,x\in\R^N,\,p\in B_{R+2}\,\},
\]
which is finite thanks to (B3). Therefore we can apply Proposition
\ref{a prop preliminary} to find a continuous $\underline
H:\T^\infty\times\R^N\to\R$ such
\[
\underline H(\tau_x\mathbf 0,p)=H_0(x,p)\qquad\hbox{for every
$(x,p)\in\R^N\times\R^N$,}
\]
where $(\tau_x)_{x\in\R^N}$ is the group of translations on
$\T^\infty$ associated via \eqref{a eq translations} to a suitably
chosen sequence $(\lambda_n)_n$ of rationally independent vectors
of $\R^N$. In particular, $\underline H$ satisfies conditions
(B2), (B3), (B4) on a dense subset of $\T^\infty\times\R^N$, hence
on the whole $\T^\infty\times\R^N$ by the continuity of
$\underline H$. The assertion readily follows with
$\Omega:=\T^\infty$ by setting
\[
H(x,p,\omega)=H(\tau_x\omega ,p)\qquad\hbox{for every
$(x,p,\omega)\in\R^N\times\R^N\times\Omega$,}
\]
and by choosing $\omega_0=\mathbf 0$.\qed
\end{section}
\end{appendix}

\end{document}